\documentclass[12pt]{article}
\usepackage{amscd}
\usepackage{mathrsfs}
\usepackage{amsfonts}
\usepackage{graphicx}
\usepackage{amsmath,amscd,stmaryrd,amssymb,array, amsfonts,amsmath,amssymb}
\usepackage{enumitem}

\setenumerate[1]{itemsep=0pt,partopsep=0pt,parsep=\parskip,topsep=5pt}
\setitemize[1]{itemsep=0pt,partopsep=0pt,parsep=\parskip,topsep=0pt}
\setdescription{itemsep=0pt,partopsep=0pt,parsep=\parskip,topsep=5pt}

\numberwithin{figure}{section}
 \numberwithin{equation}{section}

\newtheorem{theorem}{Theorem}[section]
\newtheorem{proposition}[theorem]{Proposition}
\newtheorem{definition}[theorem]{Definition}

\newtheorem{lemma}[theorem]{Lemma}
\newtheorem{remark}[theorem]{Remark}

\newcommand{\cC}{{\mathcal C}}
\newcommand{\cD}{{\mathcal D}}
\newcommand{\cF}{{\mathcal F}}
\newcommand{\cH}{{\mathcal H}}
\newcommand{\cG}{{\mathcal G}}
\newcommand{\cO}{{\mathcal O}}

\newcommand{\cM}{{\mathcal M}}
\newcommand{\cN}{{\mathcal N}}
\newcommand{\cU}{{\mathcal U}}
\newcommand{\cV}{{\mathcal V}}

\newcommand{\cS}{{\mathcal S}}
\newcommand{\cT}{{\mathcal T}}
\newcommand{\cW}{{\mathcal W}}

\newcommand{\cK}{{\mathcal K}}

\newcommand{\sC}{{\mathscr C}}

\newcommand{\sP}{{\mathscr P}}

\newcommand{\sK}{{\mathscr K}}

\newcommand{\sX}{{\mathscr X}}

\newcommand{\mB}{\mb{B}}

\def\be{\begin{equation}}
\def\ee{\end{equation}}
\def\ba{\begin{array}}
\def\ea{\end{array}}
\def\benu{\begin{enumerate}}
\def\eenu{\end{enumerate}}
\def\bt{\begin{theorem}}
\def\et{\end{theorem}}
\def\bp{\begin{proposition}}
\def\ep{\end{proposition}}
\def\bl{\begin{lemma}}
\def\el{\end{lemma}}
\def\br{\begin{remark}}
\def\er{\end{remark}}

\def\b{\beta}
\def\De{\Delta}
\def\de{\delta}
\def\pa{\partial}

\def\lam{\lambda}
\def\Lam{\Lambda}

\def\Upsi{\Upsilon}

\def\ve{\varepsilon}
\def\sig{\sigma}

\def\gam{\gamma}

\def\a{\alpha}

\def\.{\cdot}

\def\R{\mathbb{R}}

\def\A{\forall}
\def\ol{\overline}

\def\Cap{\bigcap}\def\Cup{\bigcup}

\def\ra{\rightarrow}

\def\~{\tilde}
\def\8{\infty}
\def\X{\times}
\def\({\left(}
\def\){\right)}

\def\E{\exists}
\def\mb{\mbox}
\def\emp{\emptyset}
\def\sm{\setminus}

\def\Hs{\hspace{1cm}}\def\hs{\hspace{0.5cm}}
\def\Vs{\vskip8pt}\def\vs{\vskip4pt}

\def\({\left(}\def\){\right)}

\begin{document}

\begin{center}
{\bf\Large Global Dynamic Bifurcation of Semiflows\\[1ex] and Nonlinear Evolution Equations\footnote{ This work was supported by the National Natural Science Foundation of China (11471240, 11771324).}}
\end{center}

\centerline{Luyan Zhou}
\begin{center}
{\footnotesize
{Center for Applied Mathematics,  Tianjin University\\
     Tianjin 300072,  China\\
 {\em E-mail}: zhouly@tju.edu.cn}}
\end{center}
\centerline{Desheng  Li\footnote{Corresponding author.}
}

\begin{center}
{\footnotesize
{School of Mathematics,  Tianjin University\\
          Tianjin 300072,  China\\

{\em E-mail}:  lidsmath@tju.edu.cn}}
\end{center}

{\footnotesize
\noindent{\bf Abstract.}  We establish new global bifurcation theorems for dynamical systems in terms of local semiflows on complete metric spaces. These theorems are applied to the nonlinear evolution equation $u_t+A u=f_\lambda(u)$ in a Banach space $X$,  where $A$ is a sectorial operator with compact resolvent. Assume that $0$ is always a trivial stationary solution of the equation.  We show that the global dynamic  bifurcation branch $\Gamma$ of a bifurcation point $(0,\lambda_0)$ either meets another bifurcation point $(0,\lambda_1)$, or is unbounded, completely extending the well-known Rabinowitz Global Bifurcation Theorem on operator equations to nonlinear evolution equations without any restrictions on the crossing number. In the case where $f_\lambda(u)=\lambda u+f(u)$, due to the {\em nonnegativity} of the Conley index we can even prove a stronger conclusion  asserting  that only one possibility occurs for $\Gamma$, that is,  $\Gamma$ is necessarily unbounded. This result can be expected to help us have  a  deeper understanding of the dynamics of nonlinear evolution equations.

As another example of applications of the abstract bifurcation theorems, we also discuss the bifurcation and the existence of nontrivial solutions of the elliptic equation $-\Delta u=f_\lambda(u)$ on a bounded domain in $\mathbb{R}^n$ ($n\geq 3$) associated with the homogenous Dirichlet boundary condition. Some new results with global features are  obtained.

 \Vs
\noindent{\bf Keywords:} Semiflow, evolution equation, global dynamic bifurcation, Conley index.

\vs \noindent{\bf 2010 MSC:} 34C23, 34K18, 37B30, 37G10, 37K50.

\vs \noindent{\bf Running Head:} Global Dynamic Bifurcation of Semiflows.


}

\newpage

\section{Introduction}

The dynamic bifurcation theory concerns the changes in qualitative or topological structures of limiting motions such as equilibria, periodic solutions, homoclinic orbits, heteroclinic orbits and invariant tori for nonlinear evolution equations as some relevant   parameters in the equations vary.
Historically, the subject can be traced back  to the very earlier work of Poincar$\acute{\mb e}$ \cite{Poin}
around 1892. It  is now a fundamental tool to study nonlinear problems which  enables us to understand how and when a system organizes new states and patterns near the original ``\,trivial\,'' one when the control parameter crosses a critical value.

A relatively simpler case of dynamic bifurcation is that of the bifurcations from equilibria.
Generally speaking,   there are two  typical kinds of such bifurcations  in the classical bifurcation  theory. One is the bifurcation from equilibria to equilibria (static bifurcation), and the other is from equilibria to periodic solutions (Hopf bifurcation). The former usually requires a ``crossing odd-multiplicity'' condition, 
 and has been extensively studied in the past decades; see e.g. Chow and Hale \cite{Chow} and Kielh$\ddot{\mb o}$fer  \cite{Kie}. A well-known classical  result in this line  is the celebrated  Rabinowitz Global Bifurcation Theorem. Situations become very complicated  if one drops  the ``crossing odd-multiplicity'' condition mentioned above. If the system under consideration is a gradient one, then using the bifurcation theory  on potential operator equations (see e.g.  Chang and Wang \cite{CW}, Kielh$\ddot{\mb o}$fer  \cite{Kie,Kie2}, Rabinowitz \cite{Rab2} and Schmitt and Wang \cite{SW}), one can still obtain local bifurcation results, whereas  global bifurcation remains an open problem.

  The Hopf bifurcation theory focuses on the case where there is exactly  a pair of conjugate eigenvalues  of the linearized equation cross the imaginary axis, and was fully developed in the 20{\em th} century.  There has been  a vast body of literature on how to determine Hopf bifurcation for nonlinear systems arising in applications. One can also find some nice results on  global Hopf bifurcation  in Alexander and York \cite{AY}, Chow and Mallet-Paret \cite{Chow2}, Fiedler \cite{Fied}, Sanjurio \cite{san3} and Wu \cite{Wu}, etc.

This paper is mainly concerned with the general case of bifurcations from equilibria, in which the {\em crossing number} of a system at a critical value of the control parameter  (the number of  eigenvalues  of the linearized equation crossing the imaginary axis) may be even and  greater than   two.
A particular but important case in such a situation  is the attractor bifurcation, which was first introduced by Ma and Wang in \cite{MW0} (see \cite{MW1} for a complete statement), and was  further developed by the authors into a dynamic transition theory in \cite{MW5}.  More abstract results concerning attractor bifurcation can also be found in Sanjurjio  \cite{san3}. Note that  the attractor bifurcation  theory fails to work if the trivial equilibrium solution $\theta$ of a system is neither an attractor nor a repeller of the system restricted on the center manifold of the equilibrium at a critical value $\lam_0$ of the control parameter.   However,  dynamic  bifurcation always occurs as long as the crossing number is nonzero; see e.g. Rybakowski \cite[Chap. II]{Ryba}, Ward \cite{Ward1} and Li et al. \cite{LW}.

The motivation of this work mainly comes from  \cite{LW}, in which  the authors  performed  a systematic study on the  dynamic  bifurcation of nonlinear evolution equations in terms of invariant-set bifurcation. In addition to a precise description on local dynamic bifurcation, they established a global dynamic bifurcation theorem (see \cite[Theorem 6.3]{LW}) for the equation \be\label{e:1.1}
u_t+A u=f_\lam(u)
\ee in a Banach space $X$ (without any restriction on the crossing number), where  $A$ is a sectorial operator in $X$ with { compact resolvent},  $f_\lam(\.)\in C^1(X^\a,X)$ for some $\a\in[0,1)$ ($X^\a$ denotes the fractional powers of $X$), and $f_\lam(0)=0$ for all $\lam\in\R$.
Let $(0,\lam_0)$ be a dynamic bifurcation point. Denote $\Gamma$ the  dynamic (invariant-set) bifurcation branch of $(0,\lam_0)$ (in the terminology of \cite{LW}) in a given neighborhood $\Omega$ of $(0,\lam_0)$ in $X^\a\X\R$.  Informally speaking,  the global bifurcation theorem in \cite{LW} states  that one of the following alternatives occurs: (1) \,$\Gamma$ is unbounded or meets the boundary $\pa\Omega$ of $\Omega$; \,(2) \,$\Gamma$ returns back to the point $(0,\lam_0)$; and (3) \,$\Gamma$ meets the trivial solution branch at another point $(0,\lam_1)$ with $\lam_1\ne\lam_0$.
 While the  theorem provides some interesting information on the dynamics of a system, compared with the Rabinowitz Global Bifurcation Theorem, we find that it has two obvious drawbacks. One is that  unlike in the case of the Rabinowitz's theorem which involves  only two alternatives, the theorem does not exclude the possibility  that $\Gamma$ may return back to the bifurcation point  $(0,\lam_0)$. The other is that in case alternative (3) occurs, the theorem  gives no information on  $\lam_1$ other than that $\lam_1\ne\lam_0$. Therefore it is not clear  how far the bifurcation branch $\Gamma$ can go.
 This is again different from that in the Rabinowitz's theorem, in which it is  known that a global static bifurcation branch necessarily  crosses at least two distinct  eigenvalues of the corresponding linear operator unless it is unbounded (or meets the boundary of a given domain).  The above mentioned  drawbacks may give a heavy discount to the theorem in applications.


In this present work we present  some new global dynamic bifurcation results by using a weaker notion of bifurcation branch different from the one used in \cite{LW}. First, we establish some abstract results  in the frame work of local semiflows on complete metric spaces. Let $\Phi_\lam$ ($\lam\in\R$) be a family of asymptotically compact local semiflows on a complete metric space $X$. Suppose that $\theta\in X$ is an equilibrium point of $\Phi_\lam$ for all $\lam\in\R$. Denote $\Upsi$ the set of dynamic bifurcation values of $\Phi_\lam$ ($\lam\in\R$).  Let $(\theta,\lam_0)$ be an essential dynamic  bifurcation point (i.e., there exists $\ve>0$ such that $h(\Phi_{\lam_0-\ve},\theta)\ne h(\Phi_{\lam_0+\ve},\theta)$, where $h(\Phi_{\lam},\theta)$ denotes the Conley index of $\theta$), and let $\Gamma$ be the global dynamic bifurcation branch of $(\theta,\lam_0)$ in $X\X\R$ (in the terminology of the present work).

 Our first global bifurcation theorem (Theorem \ref{t:3.2}) asserts that  if  $\Gamma$  is bounded then it meets $\{\theta\}\X\R$ at another point $(\theta,\mu_0)$ with $\mu_0\ne\lam_0$; furthermore,  there is at least
a bifurcation value $\lam_1\ne\lam_0$ between $\lam_0$ and $\mu_0$.   This result  can be seen as a generalization of an earlier one of Ward; see \cite[Theorem 4]{Ward1}. If we further assume that each bifurcation value is isolated in $\Upsi$; moreover, for any  $\lam\notin \Upsi$, we have
\be\label{e:1.2}
h(\Phi_\lam,\theta)=\Sigma^p
\ee
for some $p\geq 0$, where $\Sigma^p$ denotes the homotopy type of the pointed sphere $(S^p,x_0)$, then it can be shown that either  $\Gamma$ is unbounded, or  it   meets another bifurcation  point $(\theta,\lam_1)$ with $\lam_1\ne\lam_0$. See Theorem \ref{t:3.3} for details.

In addition to the above hypotheses,  let us now  assume that  the Conley index along the trivial solution branch satisfies a stronger  condition:
For any bounded interval  $[a,b]\subset\R$ with $a,b\notin\Upsi$, we have
 \be\label{e:1.3}
 h(\Phi_a,\theta)\ne h(\Phi_b,\theta)
 \ee
 as long as $(a,b)\cap\Upsi\ne\emp$. Then we can  prove that  the dynamic  bifurcation branch $\Gamma$ is necessarily unbounded (see Theorem \ref{t:3.7}).

At first glance,  conditions \eqref{e:1.2} and  \eqref{e:1.3} seem to be quite restrictive.  However, as we will see in Section 5, they can be  naturally satisfied by a  nonlinear evolution equation as in  \eqref{e:1.1}.
As a direct application of the above  theorems, we immediately
conclude  that under some reasonable additional assumptions, the global dynamic  bifurcation branch $\Gamma$ of a bifurcation point $(0,\lam_0)$ of   \eqref{e:1.1} either meets another bifurcation point $(0,\lam_1)$, or is unbounded in $X^\a\X\R$ (Theorem \ref{gbt1}). This result can be seen as a dynamical version of the Rabinowitz Global Bifurcation Theorem on operator equations. It is worth mentioning that in our case, we need not require  the crossing number to be odd.

If \eqref{e:1.1} takes a slightly simpler form, say,
\be\label{e:1.4}
u_t+A u=\lam u+f(u),
\ee
then  the Conley index is always {\em increasing} along the trivial solution branch. As we have mentioned above,  in such a case condition \eqref{e:1.3} is successfully fulfilled. By virtue of our abstract global bifurcation theorems, we conclude  that there is  only one possibility left  for $\Gamma$, that is,  $\Gamma$ is necessarily unbounded in $X^\a\X\R$ (see Theorem \ref{gbt2}). This result is somewhat different from that in the situation  of the Rabinowitz's theorem,  and  may help us have a  deeper understanding of the dynamics of  evolution equations.

Finally, as another example of applications of our abstract bifurcation results, we consider  the elliptic equation
\be\label{e:1.5}-\De u=f_\lam(u)\ee on a bounded domain in $\R^n$ ($n\geq 3$) associated with Dirichlet boundary condition.
Such problems aroused much interest in the past decades. In case $f_\lam(s)=\lam s+o(|s|)$ (as $|s|\ra 0$),  if $\lam_0$ is an eigenvalue of the operator $A=-\De$  of odd multiplicity, the Rabinowitz's theorem  enables  us to obtain  some global  bifurcation results.
However, if $\lam_0$ is of  even multiplicity, then the Rabinowitz's theorem fails to work. In such a case the bifurcation theory on potential operators (see e.g. \cite{CW,Kie,Kie2,Rab2} and \cite{SW}) allows us to take a step, but in general only some local  bifurcation results can be derived. A typical result in this line is as follows: Suppose  the trivial solution is  isolated at $\lam=\lam_0$. Then either  there is  a one-sided neighborhood $\Lam_1$ of $\lam_0$ such that for each $\lam\in\Lam_1\sm\{\lam_0\}$,  the problem has at least two distinct nontrivial solutions, or there is  a two-sided neighborhood $\Lam_2$ of $\lam_0$ such that for each $\lam\in\Lam_2\sm\{\lam_0\}$, the problem has at least one nontrivial solution; see e.g. \cite{CW, Rab2}.

In this present work we will try to  exhibit some global features of bifurcation for  such problems at any eigenvalue of either odd or even  multiplicity. This is summarized in Theorem \ref{t:4.1}.
For instance,  consider the  equation
\be\label{e:1.6}
-\Delta u=\lam u+\a|u|^{p-1}u+ \beta |u|^{q-1}u,
\ee
where $1<q<p<(n+2)/(n-2)$,  and $\a,\b\in\R$ are constants with $\a\ne 0$ ($\a,\b$ may be either positive or negative). As a particular case of Theorem \ref{t:4.1} we have the following interesting   result.

\bp\label{t:1.1} For  each eigenvalue $\mu_k$ of $A=-\De$, there is an interval $\Lambda$ containing  $\mu_k$  such that  $(\ref{e:1.6})$ has at least a  nontrivial solution $u_\lam\ne 0$ for all $\lam \in \Lambda\setminus\sigma(A)$. Moreover,  one of the following alternatives  occurs:
  \begin{enumerate}
\item[$(1)$] There is a bounded sequence $\lam_m\in \Lambda$ such that  $||u_{\lam_m}||\rightarrow \infty$ as $m\rightarrow\8$.
\vs
\item[$(2)$]  $\Lambda$ contains either  the interval $(-\8,\mu_{k}]$ or the interval  $[\mu_k,\8)$.
\end{enumerate}
\ep

 This paper is organized as follows. In Section 2 we make some  preliminaries, and in Section 3 we introduce the notion of a dynamic  bifurcation branch and discuss basic properties of bifurcation branches. In Section 4 we establish  global dynamic bifurcation theorems for local semiflows on metric spaces. Section 5 is devoted to the global dynamic bifurcation of the evolution equation \eqref{e:1.1}, and  Section 6  consists of some argument on global features of bifurcation of the  elliptic equation \eqref{e:1.5}. In the Appendix part we present a result on perturbation of a sectorial operator with compact resolvent.

\section{Preliminaries}

This section is concerned with some preliminaries.
\subsection{Basic topological notions and facts}
\Vs

Let $X$ be  a complete  metric space with metric $d(\cdot,\cdot)$.

 Let $A$ and $B$ be nonempty subsets of $X$.
The {\em distance} $d(A,B)$ between $A$ and $B$ is defined as
 $$
 d(A,B)=\inf\{d(x,y)|\,\,\,x\in A,\,\,y\in B\},
 $$
 and the {\em Hausdorff semi-distance} and {\em Hausdorff distance} of $A$ and $B$ are defined, respectively,
as
$$d_{\mbox{\tiny H}}(A,B)=\sup_{x \in A}d(x,B),\hs
\delta_{\mbox{\tiny H}}(A,B) = \max\left\{d_{\mbox{\tiny H}}(A,B),
d_{\mbox{\tiny H}}(B,A)\right\}.
$$
 We also assign
$d_{\tiny\mb{H}}(\emp,B)=0$.

 The
closure, interior and boundary of $A$ are denoted, respectively, by $\ol A$, int$\,A$ and $\pa A$. A
subset $U$ of $X$ is called a {\em neighborhood} of $A$, if $\ol A\subset \mbox{int}\,U$. The {\em $\ve$-neighborhood} $\mB(A,\ve)$ of $A$ is defined to be the set $\{y\in
X|\,\,d(y,A)<\ve\}.$
\vs


\bl\label{l:2.2}\cite{Rab} Let $X$ be a compact metric space, and let $A$ and $B$ be two disjoint closed subsets of $X$. Then either there exists a subcontinuum $C$ of $X$ such that
$$\ba{ll}
A\cap C\ne \emp\ne B\cap C,\ea
$$
or $X=X_A\cup X_B$, where $X_A$ and $X_B$ are disjoint compact subsets of $X$ containing $A$ and $B$, respectively.

\el

\bl\label{l:2.3} $($\cite{CV}, pp. 41$)$\, Let $X$ be a compact metric space. Denote $\sK(X)$ the family of compact subsets of $X$ which is equipped with the Hausdorff metric $\de_{\mbox{\tiny H}}(\.,\.)$. Then $\sK(X)$ is a compact metric space.
\el

\subsection{Wedge product of pointed spaces}
\Vs

Let $(X,x_0)$ and $(Y,y_0)$ be two pointed  spaces. The {\em wedge product} $(X,x_0)\vee (Y,y_0)$
is defined as $$(X,x_0)\vee (Y,y_0)=\(W,\,(x_0,y_0)\),
\hs\mb{where $W=X\X\{y_0\}\cup \{x_0\}\X Y$}.$$
Denote  $[(X,x_0)]$ the {\em homotopy type} of a  pointed space $(X,x_0)$.
Since the operation ``$\vee$\,'' preserves homotopy equivalence relations, it can be naturally extended to the homotopy types of pointed spaces. Specifically, we have
$$[(X,x_0)]\vee [(Y,y_0)]=\left[(X,x_0)\vee(Y,y_0)\right].$$ It is a basic knowledge that $[(X,x_0)]\vee \ol0=[(X,x_0)]$ for any pointed space $(X,x_0)$,  where $\ol0$ denotes the homotopy type of the one-point space $(\{p\},p)$.
See e.g. Hatcher \cite{Hat} for details of this part.

Denote $\Sigma^m$ $(m\geq0)$ the homotopy type of the  pointed $m$-dimensional sphere.
\bl\label{l£»2.2-1}\cite[Chap. I, Lemma 11.7]{Ryba}
\,Let $(X,x_0)$ and $(Y,y_0)$ be two pointed spaces. 
If $[(X,x_0)]\vee [(Y,y_0)]=\Sigma^m$ for some $m\geq0$, then either $[(X,x_0)]=\ol0$ or $[(Y,y_0)]=\ol0$.
\el

\subsection{Local semiflow and basic dynamical concepts}
\Vs

Let $X$ be a complete metric space with metric $d(\.,\.)$.
 A  {\em local semiflow} $\Phi$ on $X$ is a
continuous mapping from an open subset $\cD_\Phi$ of $\R^+\X X$  to $X$  satisfying that (1)  \,for any $x\in X$, there is a number $T_x\in(0,\8]$  such that $$(t,x)\in \cD_\Phi \Longleftrightarrow t\in[0,T_x)\,;$$
and (2) \,$\Phi(0,\.)=id_X$, and
$$
\Phi(s+t,x)=\Phi\(t,\Phi(s,x)\), \Hs x\in X,\, \,s,t\geq0,\,\,(s+t,x)\in\cD_\Phi.
$$
The number $T_x$ in (1)  is called  the {\em escape time} of $\Phi(t,x)$.
\vs
 Let $\Phi$ be a given local semiflow on $X$. For convenience,  we will write $$\Phi(t,x)=\Phi(t)x.$$

Let $M$ be a subset of $X$. Given $t> 0$, denote
$$
\Phi(t)M=\{\Phi(t)x|\,\,\,x\in M,\,\,t<T_x\}.
$$
We also write $\Cup_{s\in[0,t]}\Phi(s)M=\Phi([0,t])M$.

\vs
We say that  {\em $\Phi$ does not explode} in $M$, if  $T_x=\8$ whenever $\Phi([0,T_x))x\subset M.$
  $M$  is said to be  {\em admissible} (see Rybakowski \cite[pp. 13]{Ryba}), if for any sequences $x_n\in M$ and $t_n\ra \8$ with $\Phi([0,t_n])x_n\subset M$ for all $n$,  the sequence $\Phi(t_n)x_n$ has a convergent subsequence.
$M$ is said to be {\em strongly admissible}, if it is admissible and moreover, $\Phi$ does not explode in $M$.

\begin{definition}\label{d2.10} $\Phi$ is said to be {asymptotically compact} on $X$, if each bounded subset $B$ of $X$ is { strongly admissible}.
\end{definition}

Since we are working in a space which may be of infinite dimensional, from now on   we always assume that
\benu\item[{\bf (AC)}] \,{\em $\Phi$ is asymptotically compact on $X$.}\eenu
\noindent This compactness  requirement is fulfilled   by a large number of examples in applications.

\vs

A set $M\subset X$ is said to be {\em invariant}, if $\Phi(t)M=M$ for all $t\geq0$.

The  proposition below collects some well-known  basic properties about  bounded invariant sets that will be frequently used in this paper.

\bp\label{r:2.5}Let $M$ be a bounded invariant set of $\Phi$. Then $(1)$ \,$M$ is precompact; and $(2)$ \,$\ol M$ is a compact invariant set of $\Phi$.
\ep


\noindent$\bullet$ ~Given $U\subset X$, we  denote
$
I(\Phi,U)
$
the {\em maximal compact invariant set} of $\Phi$ in $U$, if any.


\br\label{r:2.6} In general $I(\Phi,U)$ may not exist. For instance, let $\Phi$ be the semiflow  generated by the linear scalar equation $x'(t)=0$. Then each subset of $\R$ is an invariant set of $\Phi$. Hence $I(\Phi,\R)$ does not exist. Note also that the whole phase space $X=\R$ is the largest invariant set $\Phi$. Thus in general  one should distinguish $I(\Phi,U)$ from the {largest invariant set}  of $\Phi$ in $U$.

However, if $U$ is a  bounded closed subset of  $X$, then by virtue of Proposition \ref{r:2.5} it is easy to see that $I(\Phi, U)$ exists and coincides with the largest invariant set of $\Phi$ in $U$ $($which may be void$)$.
\er

\vs A {\em trajectory} on an interval $J$  is a continuous  mapping $\gamma: J\ra X$ such that
$$
\gamma(t)=\Phi(t-s)\gamma(s),\Hs \A \,t,\, s\in J,\,\,\,t\geq s.
$$
The set $\mb{orb}(\gamma)=\{\gamma(t)|\,\,t\in J\}$ is called the {\em orbit} of $\gamma$

A {trajectory} $\gamma: \R\ra X$ is simply called a {\em complete trajectory}.
The {\em $\omega$-limit set} $\omega(\gamma)$ and {\em $\alpha$-limit set} $\alpha(\gamma)$ of    a { complete trajectory} $\gamma$ are defined  as
$$
\omega(\gamma)=\{y\in X|\,\,\,\mb{$\E$ $t_n\ra\8$ such that $\gamma(t_n)\ra y$}\},
$$
$$
\a(\gamma)=\{y\in X|\,\,\,\mb{$\E$  $t_n\ra-\8$ such that $\gamma(t_n)\ra y$}\}
$$
respectively.

\subsection{Conley index}
\Vs

Let us  recall  briefly some basic notions and  results in the  Conley index theory. The interested reader is referred to \cite{Conley,Mis} and \cite{Ryba} for details.


A compact invariant set $S$ of $\Phi$ is said to be {\em isolated}, if there exists a neighborhood $N$ of $S$ such that $S=I(\Phi,\ol N)$, namely,
$S$ is the  maximal compact invariant set in $\ol N$. Consequently $N$ is called
 an {\em isolating neighborhood} of $S$.

Let $S$ be  an isolated compact invariant set.  A pair of bounded closed subsets $(N,E)$ is called an {\em  index pair} of $S$, if (1)\, $N\sm E$ is an isolating neighborhood of $S$;
(2) \,$E$ is $N$-invariant, namely, for any $x\in E$ and $t\geq 0$, $$\Phi([0,t])x\subset
N\Longrightarrow\Phi([0,t])x\subset E;$$
(3)\, $E$ is an exit set of $N$. That is,  for any $x\in N$, if $\Phi(t_1)x\not\in N$ for some
$t_1>0$, then there exists  $0\leq t_0\leq t_1$ such that
    $\Phi(t_0)x\in E.$

\br 
Index pairs in the terminology of \cite{Ryba} need not be bounded. However, the bounded ones are sufficient for our purposes here.
\er

\begin{definition}  
The {homotopy Conley index} of isolated compact invariant set $S$ of $\Phi$, denoted by $h(\Phi,S)$,  is defined to be the homotopy
type $[(N/E,[E])]$ of the pointed space $(N/E,[E])$ for any index pair $(N,E)$ of $S$.
\end{definition}
\br
We assign  $h(\Phi,\emptyset)=\ol0$.
\er
\bl\label{l:2.4.1} \cite{Ryba}
Let $S_1$, $S_2$ be two isolated compact invariant sets of $\Phi$ with $S_1\cap S_2=\emptyset$. Then
$$h(\Phi, S_1\cup S_2)=h(\Phi, S_1)\vee h(\Phi, S_2).$$
\el


 Let $\{\Phi_\lam\}_{\lam\in\Lam}$ be a family of semiflows, where $\Lam$ is a complete  metric space. Assume $\Phi_\lam(t)x$ is continuous in $(t,x,\lam)$.
Denote $\~\Phi$ the {\em skew-product  flow} of the family $\{\Phi_\lam\}_{\lam\in\Lam}$ on $\sX=X\X \Lam$ defined as follows:
\begin{equation}\label{spf}
\~\Phi(t)(x,\lam)=\(\Phi_\lam(t)x,\lam\),\Hs (x,\lam)\in X\X \Lam.
\end{equation}

\noindent $\bullet$~~~Let    $\cF\subset \sX$.  For any $\lam\in\Lam$,  denote $\cF[\lam]$ the {\em $\lam$-section}  of $\cF$,
$$
\cF[\lam]=\{x|\,\,(x,\lam)\in \cF\}.
$$

The following continuation result is a slightly modified version of Ward \cite[Theorem 2]{Ward1}, which seems to be more natural and convenient in applications.

\bt\label{t:2.14} Let $\Lam\subset\R$ be a compact interval. Suppose that $\~\Phi$ is asymptotically compact. Let $\cS$ be  an  isolated compact invariant set of $\~\Phi$.
Then $$h(\Phi_\lam,\cS[\lam])\equiv \mb{const.},\Hs \lam\in\Lam.$$
\et

\noindent{\bf Proof.} We give a self-contained proof for  the reader's convenience, which is simpler than that of \cite[Theorem 2]{Ward1}.

 Using the compactness of $\cS$ one can easily verify that the $\lam$-section $\cS[\lam]:=S_\lam$ is upper semicontinuous in $\lam$, namely, $d_{\mbox{\tiny H}}(S_{\lam},S_{\lam_0})\ra0$ as $\lam\ra\lam_0\in\Lam$.

Take a bounded closed isolating neighborhood $\cN$ of $\cS$ in $X\X\Lam$. Then if $S_{\lam}\ne\emp$,  $\cN[\lam]:=N_\lam$  is an isolating neighborhood  of $S_\lam$.

Let $\lam_0\in \Lam$. If $S_{\lam_0}\ne\emp$, then by the stability of isolating neighborhoods we deduce that  $N_{\lam_0}$ is  an isolating neighborhood  of $\Phi_{\lam}$ for all $\lam$ in a small neighborhood $U_\ve=(\lam_0-\ve,\lam_0+\ve)\cap\Lam$ of $\lam_0$ in $\Lam$. We show  that \be\label{e:2.a}I(\Phi_\lam,N_{\lam_0})=S_\lam,\Hs \lam\in U_\ve,\ee provided $\ve>0$ is sufficiently small.

Indeed, let $I(\Phi_{\lam},N_{\lam_0})=\~S_\lam$. Then by a very standard argument we can verify  that $d_{\mbox{\tiny H}}(\~S_\lam,S_{\lam_0})\ra0$ as $\lam\ra\lam_0$. This implies that
$$
d_{\mbox{\tiny H}}\(\~S_\lam\X\{\lam\},\,S_{\lam_0}\X\{\lam_0\}\)\ra0\hs \mb{as $\lam\ra\lam_0$}.$$ Thus there exists $\ve>0$ such that  $\~S_\lam\X\{\lam\}\subset \cN$ for $\lam\in U_\ve$. Therefore $\~S_\lam\subset N_\lam$. As $N_\lam$ is an isolating neighborhood of $S_\lam$, it follows that
$\~S_\lam\subset S_\lam$.
On the other hand, since $d_{\mbox{\tiny H}}(S_{\lam},S_{\lam_0})\ra0$ as $\lam\ra\lam_0$, we have $S_\lam\subset N_{\lam_0}$ for $\lam\in U_\ve$ provided that  $\ve$ is sufficiently small. Consequently
$$
S_\lam \subset I(\Phi_\lam, N_{\lam_0})=\~S_\lam.
$$
Hence we see that \eqref{e:2.a} holds true.

By virtue of    \cite[Chap.\,\,1,  Theorem 12.2.]{Ryba} and \eqref{e:2.a}  we deduce that $$h(\Phi_\lam,S_\lam)\equiv\mb{const.},\Hs\lam\in U_\ve.$$

Now assume  $S_{\lam_0}=\emp$. Then by the upper semicontinuity of $S_\lam$ it is trivial to check  that there is a small neighborhood $U_\ve$ of $\lam_0$  such that $S_\lam=\emp$ for all $\lam\in U_\ve$. Hence $h(\Phi_\lam,S_\lam)\equiv\ol0$ on $U_\ve$.

 In conclusion, for each  $\lam_0\in\Lam$, one can always  find  a neighborhood $U_\ve$ of $\lam_0$ such that $h(\Phi_\lam,S_\lam)\equiv\mb{const.}$ on $U_\ve$.

 Fix a $\lam^*\in\Lam$, and set
 $$
 \Lam_0=\{\lam\in\Lam|\,\,\,h(\Phi_\lam,S_\lam)=h(\Phi_{\lam^*},S_{\lam^*})\}.
 $$
 Using what we have proved above, it is trivial to check that $\Lam_0$ is both  open and  closed in $\Lam$. The connectedness of $\Lam$ then asserts that  $\Lam_0=\Lam$. \,$\Box$

\br\label{r2.1}We emphasize that in Theorem \ref{t:2.14}, if $\cS[\lam_0]=\emp$ for some $\lam_0\in\Lam$ then
$h(\Phi_\lam,\cS[\lam])\equiv \ol0$ on $\Lam$.
\er

\section{Dynamic Bifurcation Branch}


Let $X$ be a complete metric space with metric $d(\.,\.)$, and let $\{\Phi_\lam\}_{\lam\in\R}$  be a family of asymptotically compact local semiflows on $X$. Assume $\Phi_\lam(t)x$ is continuous in $(t,x,\lam)$.

Set $\sX=X\X\R$,  which is equipped with the metric $\rho(\.,\.)$ defined by
\be\label{e:3.0}
\rho\((x,\lam),\,(y,\lam')\)=d(x,y)+|\lam-\lam'|, \Hs  (x,\lam),\, (y,\lam')\in\sX.
\ee
For    $\cF\subset \sX$, denote $\cF[\lam]$ ($\lam\in\R$) the {\em $\lam$-section} of $\cF$, $\cF[\lam]=\{x|\,\,(x,\lam)\in \cF\}.$

Let $\~\Phi$ be the skew-product flow of $\{\Phi_\lam\}_{\lam\in\R}$   on $\sX$.  From now on we always assume that  $\~\Phi$  is asymptotically compact.

Suppose that $\theta\in X$ is an equilibrium of $\Phi_\lam$ for all $\lam$.
For notational simplicity, we usually write $$h(\Phi_\lam,\{\theta\})=h(\Phi_\lam,\theta)$$ in case $S_0=\{\theta\}$ is an isolated invariant set of $\Phi_\lam$.


\begin{definition}\label{d:3.1}Let  $\lam_0\in\R$.  If  for any neighborhood $U$ of $\theta$ and $\ve>0$, there exists $\lam\in (\lam_0-\ve,\lam_0+\ve)$ such that $\Phi_{\lam}$ has a nonempty compact invariant set $K_\lam\subset U$ with $K_\lam\ne\{\theta\}$, then we call $\lam_0$ a $(dynamic)$ {bifurcation value} of $\{\Phi_\lam\}_{\lam\in\R}$  $($along the trivial equilibrium branch $\{\theta\}\X\R)$. Accordingly,
  $(\theta,\lam_0)$ is called a $(dynamic)$ {bifurcation point}.
\end{definition}

\noindent $\bullet$ \,We denote $\Upsilon$ the set of  bifurcation values of $\{\Phi_\lam\}_{\lam\in\R}$.

\bp $\Upsilon$ is a closed subset of $\R$.
\ep
\noindent{\bf Proof.} This is a simple consequence of the definition of bifurcation values.
\,$\Box$

\Vs If $\lam_0\in \Upsilon$ is an {\em isolated bifurcation value}, then there is a number  $\ve>0$ such that $S_0=\{\theta\}$ is an isolated invariant set of $\Phi_\lam$ for each $\lam\in (\lam_0-\ve,\lam_0+\ve)\sm\{\lam_0\}$. Thus by the continuation of the Conley index we find that
$$
h(\Phi_\lam,\theta)\equiv \mb{const}.:=h(\Phi_{\lam_0^-},\theta),\Hs \lam\in (\lam_0-\ve,\lam_0),$$
$$
 h(\Phi_\lam,\theta)\equiv\mb{ const}.:=h(\Phi_{\lam_0^+},\theta),\Hs \lam\in (\lam_0,\lam_0+\ve).
$$
\begin{definition}Let $\lam_0\in \Upsilon$ be an {\em isolated bifurcation value}. We  call $\lam_0$ an {essential  bifurcation value} of $\{\Phi_\lam\}_{\lam\in\R}$ if   $h(\Phi_{\lam_0^-},\theta)\ne h(\Phi_{\lam_0^+},\theta)$. Accordingly, $(\theta,\lam_0)$ is called an { essential  bifurcation point}.
\end{definition}


Denotes $\sC(\Phi_\lam)$ the family of {\em connected} compact invariant sets $C$ of $\Phi_\lam$ with $C\ne\{\theta\}$.
Given $\cN\subset \sX$, set  \be\label{e:sP}\ba{ll}
\sP(\cN)=\ol{\Cup \{C\X\{\lam\}\subset \cN|\,\,\,\, C\in\sC(\Phi_\lam),\,\,\lam\in\R\}}.\ea
\ee
\br\label{r:3.3} We infer from the construction of $\sP(\cN)$ that for any $(x^*,\lam^*)\in \sP(\cN)$, there is a sequence $\lam_n\rightarrow\lam^*$ such that for each $n$, $\Phi_{\lam_n}$ has a connected compact invariant set $C_n\neq \{\theta\}$ with $C_n\X\{\lam_n\}\subset \cN$ and
$\lim_{n\rightarrow \infty}d(x^*,C_n)=0.$
\er
\bp\label{p:3.3}If $\cN$ is bounded,  $\sP(\cN)$ is a compact invariant set of $\~\Phi$.
\ep
\noindent{\bf Proof.} It is easy to see that $\cK=\Cup \{C\X\{\lam\}\subset \cN|\,\,\,\, C\in\sC(\Phi_\lam),\,\,\lam\in\R\}$ is
 an  invariant set of $\~\Phi$. Thus the conclusion immediately follows from Proposition \ref{r:2.5} and the asymptotic compactness of $\~\Phi$. \,$\Box$
\Vs
\noindent $\bullet$~~Let $(\theta,\lam_0)$ be a bifurcation point. Given $\cN\subset\sX$ with $(\theta,\lam_0)\in\cN$, we denote $\Gamma_{\cN}(\theta,\lam_0)$ the $(connected)$ component of $\sP(\cN)$ containing $(\theta,\lam_0)$.

\br\label{r:3.6}If $\cN$ is bounded, then by Proposition \ref{p:3.3} we deduce that  $\Gamma_\cN(\theta,\lam_0)$ is a compact invariant set of $\~\Phi$.

\er

\begin{definition}\label{d3.1} Let $(\theta,\lam_0)$ be a bifurcation point.
  The {global dynamic bifurcation branch}   of $(\theta,\lam_0)$, denoted by $\Gamma(\theta,\lam_0)$,  is defined  as
$$\ba{ll}\Gamma(\theta,\lam_0)=\Cup_{n\geq 1}\Gamma_{\cN_n}(\theta,\lam_0),\ea
$$
where $\cN_n=\ol\mB((\theta,\lam_0),n)$ is  the ball in $\sX$ centered at $(\theta,\lam_0)$ with radius $r=n$.
\end{definition}

\br One may simply define the global dynamic  bifurcation branch of  $(\theta,\lam_0)$  to be  the component $\Gamma_\sX(\theta,\lam_0)$ of $\sP(\sX)$ containing  $(\theta,\lam_0)$. Clearly $\Gamma(\theta,\lam_0)\subset \Gamma_\sX(\theta,\lam_0)$, where $\Gamma(\theta,\lam_0)$ is given by  Definition \ref{d3.1}.

In general, $\Gamma_\sX(\theta,\lam_0)$ may be larger than  $\Gamma(\theta,\lam_0)$, and  we do not know whether the assertion $(2)$ in Proposition \ref{p:3.2} below remains valid for $\Gamma_\sX(\theta,\lam_0)$.
\er

\br\label{r:3.4}
We infer from Remark \ref{r:3.6} and  Definition \ref{d3.1} that the global bifurcation branch  $\Gamma=\Gamma(\theta,\lam_0)$ is the union of at most countably infinitely many  compact invariant sets  of $\~\Phi$.  Consequently for each $\lam$, the section $\Gamma[\lam]$ can be expressed as the union of a family of  compact invariant sets  of $\Phi_\lam$. As a result, we see that $\Gamma[\lam]$  consists of orbits of bounded complete trajectories.
\er

\bp\label{p:3.2} Let $(\theta,\lam_0)$ be a bifurcation point, and  let $\Gamma=\Gamma(\theta,\lam_0)$. Then $(1)$\, $\Gamma$ is connected; and
$(2)$ \,$\lam^*\in  \Upsilon$ whenever $\Gamma[\lam^*]=\{\theta\}$.
\ep

\noindent{\bf Proof.}  Let $\cN_n$ and  $\Gamma_{\cN_n}(\theta,\lam_0):=\Gamma_n$  be  as in Definition \ref{d3.1}.
Since $\Gamma_n$ is connected and that $(\theta,\lam_0)\in\Gamma_n$ for all  $n$, Clearly $\Gamma=\Cup_{n\geq 1}\Gamma_n$ is connected as well. Hence the assertion (1) holds true.

\vs

Now assume $\Gamma[\lam^*]=\{\theta\}$.
As $(\theta,\lam^*)\in\Gamma=\Cup_{n\geq 1}\Gamma_n$, there is a number  $m\geq 1$ such that $
(\theta,\lam^*)\in \Gamma_{m}.
$
As $\Gamma_m\subset\Gamma$, we necessarily have  \be\label{e3.00}\Gamma_m[\lam^*]=\{\theta\}.\ee

 By Remark \ref{r:3.3} there is a sequence $\lam_n\rightarrow\lam^*$ such that for each $n$, $\Phi_{\lam_n}$ has a connected compact invariant set $C_n\neq \{\theta\}$ with $C_n\X\{\lam_n\}\subset \cN_m$, such that
\be\label{e3.01}\lim_{n\rightarrow \infty}d(\theta,C_n)=0.\ee
In what follows we show that \be\label{e3.23}\lim_{n\ra\8}d_{\mbox{\tiny H}}\(C_n,\{\theta\}\)=0.\ee It then follows that $(\theta,\lam^*)$ is a bifurcation point, hence $\lam^*\in \Upsilon$, which completes the proof of the assertion (2).

We argue by contradiction and suppose \eqref{e3.23} was false. Then there would exist  a  closed neighborhood $N$ of $\theta$ 
and a subsequence of $C_n$, still denoted by $C_n$,  such that
\be\label{e3.02}C_n\backslash N\neq \emp\ee
 for each $n$. Since $\cN_m$ is closed and  bounded, by Remark \ref{r:2.6} we deduce that $\cM=I(\~\Phi,\cN_m)$ is compact. It follows that the union of all the sections $\cM[\lam]$ of $\cM$, denoted by $M$,  is precompact. As $C_n\subset M$ for all $n$, Lemma \ref{l:2.3} asserts that, up to a subsequence, $C_n$ converges in the sense of Hausdorff distance to a compact set $C$. It is trivial to check that $C$ is a connected invariant set of $\Phi_{\lam^*}$.

By \eqref{e3.01} we see that $\theta\in C$. Because $C_n\X\{\lam_n\}\subset \cN_m$, we necessarily have $C\X\{\lam^*\}\subset \cN_m$. Hence by the connectedness of $C$ and the fact that $(\theta,\lam^*)\in\Gamma_m\cap\(C\X\{\lam^*\}\)$, we deduce  that $\Gamma_m\cup(C\X\{\lam^*\})$ is a continuum in $\sP(\cN_m)$ containing $(\theta,\lam_0)$. Therefore $C\X\{\lam^*\}\subset \Gamma_m$. Thus by \eqref{e3.00} one concludes that $C=\{\theta\}$.

 On the other hand, by \eqref{e3.02} one can  easily verify  that $C\ne\{\theta\}$. This leads to a contradiction.
  \, $\Box$

  \section{Global Dynamic Bifurcation}
In this section we state and prove our abstract global bifurcation results in the framework  of local semiflows on complete metric spaces.

 Let $X$, $\{\Phi_\lam\}_{\lam\in\R}$ and $\~\Phi$ be the same as in Section 3,
   and let $\sX=X\X\R$, which is equipped with the metric $\rho(\.,\.)$ given by \eqref{e:3.0}.

 Suppose that $\theta\in X$ is an equilibrium point of $\Phi_\lam$ for all $\lam$. Let $\Upsi$ be the set of bifurcation values of $\{\Phi_\lam\}_{\lam\in\R}$ along the trivial equilibrium branch $\{\theta\}\X\R.$

\subsection{A first global dynamic  bifurcation theorem}
\Vs

 Our first result is summarized in the following theorem.

 \bt\label{t:3.2} $($A first global dynamic  bifurcation theorem$)$ 
 Let $(\theta,\lam_0)$ be an essential bifurcation point of $\{\Phi_\lam\}_{\lam\in\R}$, and let $\Gamma=\Gamma(\theta,\lam_0)$ be the global dynamic bifurcation branch. Then one of the following alternatives occurs.
 \benu\item[$(1)$] $\Gamma$  is unbounded in $\sX$;
 \item[$(2)$] $\Gamma$ meets $\{\theta\}\X\R$ at  another  point  $(\theta,\mu_0)$ with $\mu_0\ne\lam_0$. Furthermore, there is at least a bifurcation value $\lam_1\in\Upsi$ between $\lam_0$ and $\mu_0$ with $\lam_1\ne\lam_0$.
\eenu
 \et

\noindent{\bf Proof.}
We assume $\Gamma$ is bounded and show that  alternative (2) occurs.

%
First, by the boundedness of  $\Gamma$ one can pick two numbers $d,R>0$ such that
\be\label{e:4.0}
\ol\mB(\Gamma,1)\subset \ol B_R\X[-d,d]:=\cC.
\ee
Here (and below) $\mB(\Gamma,r)$ and $B_r:=\mB(\theta,r)$ denote the $r$-neighborhood of $\Gamma$ in $\sX$ and the ball in $X$ centered at $\theta$ with radius $r$, respectively. Then $\Gamma$ coincides with the component $\Gamma_\cC(\theta,\lam_0)$ of $\sP(\cC)$ containing $(\theta,\lam_0)$. By Remark \ref{r:3.6} we conclude that $\Gamma$ is a compact invariant set of the skew-product flow $\~\Phi$.

 Set
\be\label{e:dJ}J=\{\lam\in\R|\,\,(\theta,\lam)\in\Gamma\}.\ee
Then $J$ is compact. Let
\be\label{e:j2}\a=\min\{\lam|\,\,\lam\in J\},\hs \b=\max\{\lam|\,\,\lam\in J\}.
\ee
To prove the assertion (2), it suffices to check that the interval $[\a,\b]$ contains another bifurcation value $\lam_1\in\Upsi$ with $\lam_1\ne \lam_0$.

We argue by contradiction and suppose the contrary.
Then
$$
[\a,\b]\cap \Upsilon=\{\lam_0\}=J\cap \Upsilon.
$$
We claim that there exists $\ve>0$ such that
\be\label{e:3.5}
[\a-\ve,\b+\ve]\cap\Upsilon=\{\lam_0\}.\ee
Indeed, if this was false, there would exist a sequence $\ve_n\ra0$ such that for each $n$, one can find a $\lam_n\in [\a-\ve_n,\b+\ve_n]\cap\Upsilon$ with $\lam_n\ne\lam_0$. As $\lam_0$ is isolated in $\Upsilon$, there is a number $\de>0$ such that $|\lam_n-\lam_0|>\de$ for all $n$. It can be assumed that $\lam_n\ra\lam^*$. Clearly $\lam^*\ne\lam_0$. On the other hand, since $\Upsilon$ is closed, we necessarily have $\lam^*\in[\a,\b]\cap\Upsilon$, which leads to a contradiction.

By \eqref{e:4.0} it can be assumed that $$\Lam:=[\a-\ve,\b+\ve]\subset [-d,d].$$

Let
\be\label{e:j3}J_\ve=\{\lam\in\R|\,\,d(\lam,J)<\ve\},\hs {{J_\ve}^c}=\R\backslash{J_\ve}.\ee
As $J_\ve^c$ is closed, by the compactness of $\Gamma$ it is trivial to check that $F_0:= \bigcup_{\lam\in{J_\ve}^{c}}\Gamma[\lam]$ is compact.
On the other hand, by the definition of ${J_\ve}^c$ it is clear   that $\theta\not\in F_0$. Thus we have
\be\label{e3.3} d(\theta,F_0):=3\de_0>0.\ee
Pick a number $r$ with $0<r<\min(\de_0,1)$ such  that $\cU=\ol\mB(\Gamma,r)\subset\cC$ and
 \be\label{e:3.4}d(\theta,\cU[\lam])\geq2\de_0,\Hs \A\, \lam\in{J_\ve}^{c}.\ee
 Set $$\cK=\sP\(\cC\)\cap\cU.$$ Since $\sP(\cC)$ is compact (see Proposition  \ref{p:3.3}),  $\cK$ is a compact subset of $\sX$.
\vs

\begin{figure}[h!]
  \centering
  \includegraphics[width=6.6cm]{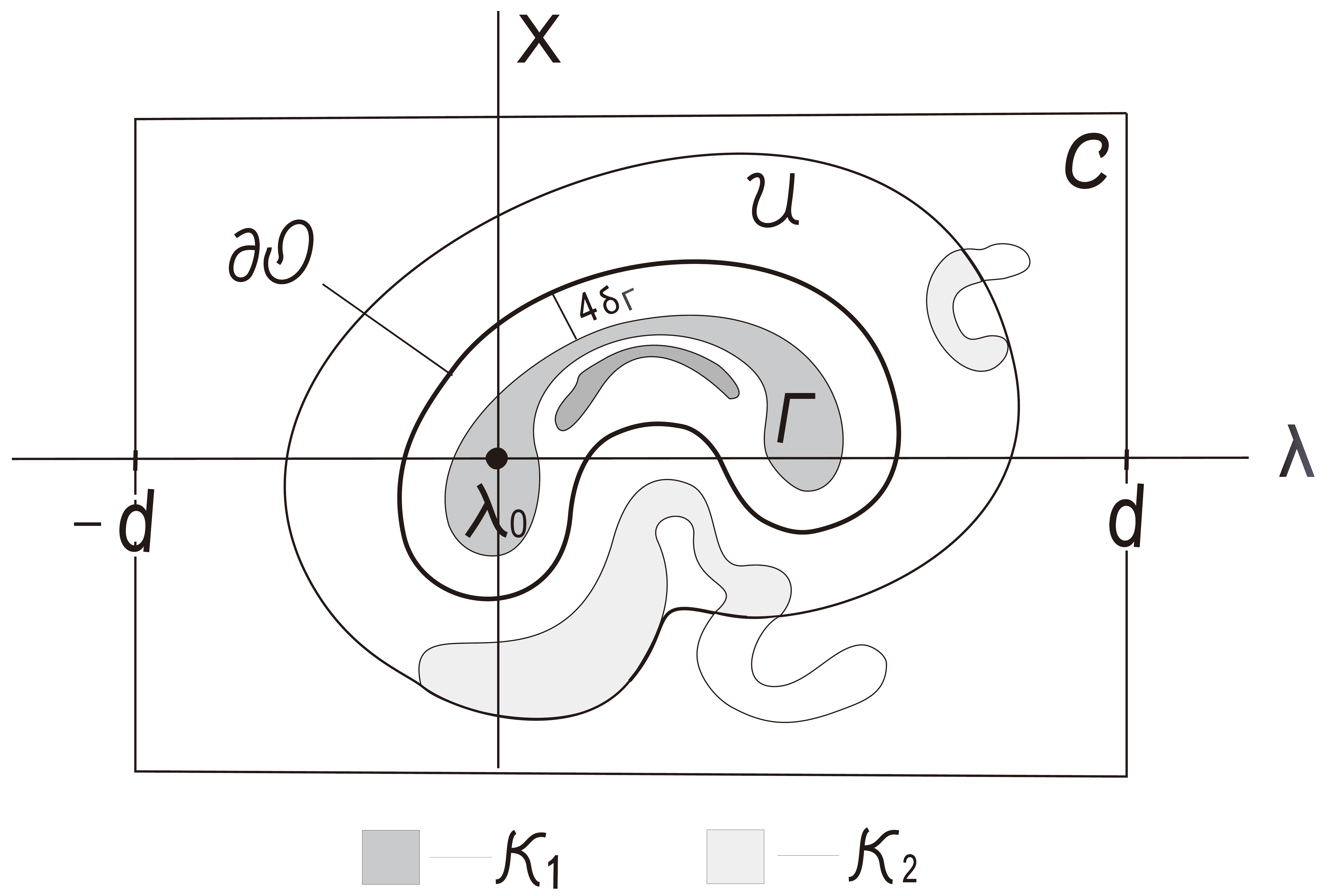}
  \caption{Separation of $\Gamma$}\label{fg3-a1}
\end{figure}

Let
$\hat{\cK}=\cK\cap \pa{\cU}.$
 Since $\pa{\cU}\cap\Gamma=\emptyset$, we have $\hat\cK\cap\Gamma=\emp$. At this point we may apply Lemma \ref{l:2.2}, the separation lemma,  to $\cK$ and its subsets $A:=\Gamma$ and $B:=\hat\cK$. Because   $\Gamma$ does not intersect any other component of $\sP(\cC)$ (and hence $\Gamma$ does not intersect any other component of $\cK$),  the first alternative in Lemma \ref{l:2.2} will not occur. Hence we deduce  that there are disjoint compact subsets $\cK_1$ and $\cK_2$ of $\cK$ with $\cK_1\cup\cK_2=\cK$ such that
$$\Gamma\subset\cK_1,\hs \hat{\cK}\subset\cK_2.$$
Clearly $\cK_1\subset \mb{int}\,\cU$; see Fig. \ref{fg3-a1}.

Take a number $\de_\Gamma$ with  $$0<\de_\Gamma<\frac{1}{8}\min\{d(\cK_1,\cK_2),\,d(\cK_1,\pa\cU)\}.$$
Set $\cO=\ol \mB(\cK_1,4\de_\Gamma)$. Then $\cO\subset \cU$. Hence by \eqref{e:3.4} we have
\be\label{e:3.8}d(\theta,\cO[\lam])>2\de_0,\Hs \A\,\lam\in {J_\ve}^{c}.\ee
Since $\cK=\cK_1\cup\cK_2$, by the choice of $\de_\Gamma$ it is  easy to see   that
$\mB(\pa\cO,2\de_\Gamma)\cap\cK=\emp.$
Consequently
\be\label{e:3.6}\ba{ll}\mB(\pa\cO,2\de_\Gamma)\Cap\sP(\cC)&=\mB(\pa\cO,2\de_\Gamma)\Cap\(\sP(\cC)\cap \cU\)\\[1ex]&=\mB(\pa\cO,2\de_\Gamma)\Cap\cK=\emp.\ea\ee
%

We claim that there exist $\sigma>0$ and $0<\mu<\ve/2$ such that
\be\label{e3.9b}B_{2\sig}\subset \cO[\lam],\Hs\A\, \lam\in J_{2\mu},\ee
where $B_r:=\mB(\theta,r)$ denotes the ball in $X$ centered at $\theta$ with radius $r$. Indeed, if the claim was not true, then for $\sigma_k=1/k$ ($k=1,2,\cdots$) there would exist  a sequence $\ve_k\rightarrow 0$ such that
$$B_{\sig_k}\varsubsetneq \cO[\lam_k]\, \,\, \mb{for some}\,\, \lam_k\in J_{\ve_k}.$$
Thus for each $k$ we can pick an $x_k\in\pa\cO[\lam_k]\cap B_{\sig_k}.$
Clearly $x_k\rightarrow \theta$ as $k\rightarrow \infty$. We may assume $\lam_k\ra \ol\lam$.
Then since $(x_k,\lam_k)\in \pa\cO$, we have
$(\theta,\ol\lam)\in  \pa{\cO}$. On the other hand, $\lam_k\in J_{\ve_k}$ implies that $\ol\lam\in J$. Hence by the definition of $J$ (see \eqref{e:dJ}) we have $$(\theta,\ol\lam)\in \Gamma\subset \mb{int}\,\cO,$$ which leads to  a contradiction\,.
\vs
In what follows we check that there exists $0<\rho<\min\{\de_0,\sigma\}$ such that
\be\label{e:3.9}I\(\Phi_\lam,\ol B_{\rho}\)=\{\theta\},\Hs \A\,\lam\in\Lam\sm {J_{2\mu}},\ee
where $\Lam=[\a-\ve,\b+\ve]$. Suppose the contrary. There would  exist sequences  $\rho_n \rightarrow 0$ and $\lam_n\in \Lam\sm {J_{2\mu}}$ such that
$I\(\Phi_{\lam_n},\ol B_{\rho_n}\)\ne\{\theta\}$ for all $n$.
It can be assumed that  $\lam_n\rightarrow \bar{\lam}\in \Lam\sm {J_{2\mu}}$.
Then $\bar{\lam}$ is a bifurcation value, which leads to a contradiction as $(\Lam\sm J_{2\mu})\cap \Upsi=\emp$. (Recall  that $\lam_0\in J\subset J_{2\mu}$ is the unique bifurcation value  in $\Lam$.)
\vs

\begin{figure}[h!]
  \centering
  \includegraphics[width=5.8cm]{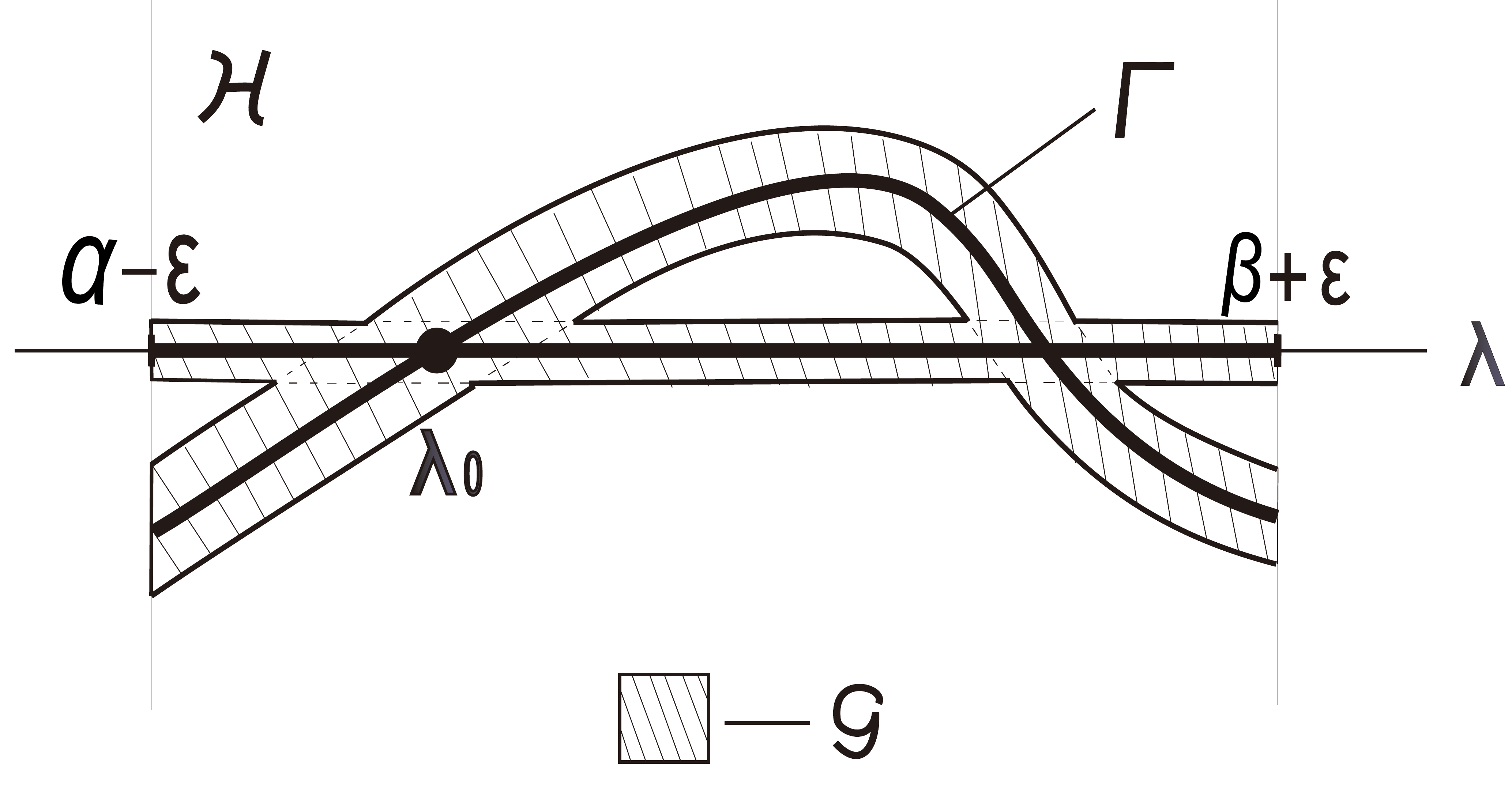}\\
  \caption{$\cG$ is an  isolating neighborhood of  $\cK$ in $\cH$}\label{fg3-2}
\end{figure}

Denote $\cF=I(\~\Phi,\cO)$. Let $\cH=X\X\Lam$, and write  $$\mb{$\cO_\cH=\cO\cap\cH$,\hs and\, $\cF_\cH=\cF\cap\cH$}.$$
Let $\cG=\cO_\cH\cup(\ol B_{\rho}\X \Lambda)$; see Fig. \ref{fg3-2}. We check  that
\be\label{e3.12b}\ba{ll}\cK:=I(\~\Phi,\cG)=\cF_\cH\cup\(\{\theta\}\X \Lambda\).\ea\ee
Indeed, it is obvious that $$\ba{ll}\cK\supset\cF_\cH\cup\(\{\theta\}\X \Lambda\).\ea$$  We show that
\be\label{e3.12c}
\cK[\lam]\subset\cO_\cH[\lam]\cup\{\theta\}= \cO[\lam]\cup\{\theta\},\Hs\lam\in\Lam,
\ee
which implies  that $\cK\subset \cF_\cH\cup\(\{\theta\}\X \Lambda\)$ and completes the proof of \eqref{e3.12b}.

Since  $\cK$ is a compact invariant set of $\~\Phi$, $\cK[\lam]$ is a compact invariant set of $\Phi_\lam$ for $\lam\in\Lam$. Thus to verify \eqref{e3.12c}, it suffices to check that for any  component $S_\lam$ of $\cK[\lam]$ ($\lam\in\Lam=[\a-\ve,\b+\ve]$), we have \be\label{e3.12d}\mb{ either $S_\lam=\{\theta\}$,\,\, or $S_\lam \subset \mb{int}\,\cO[\lam]$.}\ee
So we assume $S_\lam\ne\{\theta\}$. Then since $S_\lam$ is connected and
$
S_\lam\X\{\lam\}\subset \cG\subset \cC,
$
by the definition of $\sP(\cC)$ (see \eqref{e:sP}) we deduce that $\cS_\lam:=S_\lam\X\{\lam\}\subset\sP(\cC)$.
\eqref{e:3.6} then asserts that $$\mb{ either $\cS_\lam\subset \mb{int}\cO$, \,\,or $\cS_\lam\cap \cO=\emp$.}$$
 We prove  that the latter case can not occur. Therefore $ \cS_\lam\subset \mb{int}\,\cO$.
 It follows that $S_\lam \subset \mb{int}\,\cO[\lam]$, which justifies \eqref{e3.12d}.

 Suppose  $\cS_\lam\cap \cO=\emp$. Then since $\rho<\sig$, we deduce that
$$
 \cS_\lam\subset \cG\sm \cO=(\ol B_\rho\X\Lam)\sm \cO\subset (\mb{by \eqref{e3.9b}})\subset \ol B_\rho\X (\Lam\sm J_{2\mu}),
 $$
which implies $\lam\in \Lam\sm J_{2\mu}$ and that $S_\lam\subset \ol B_\rho$. This contradicts \eqref{e:3.9} as  $S_\lam\ne\{\theta\}$.

  We claim that $\cG$ is an isolating neighborhood of $\cK$ with respect to $\~\Phi$ restricted on $\cH$.
  Indeed, it is trivial to verify  that
\be\label{e3.12g}\pa_\cH\cG\subset\pa\cO\cup\((\pa B_{\rho}\X\Lam)\backslash{\cO}\),\ee
where $\pa_\cH\cG$ denotes the boundary of $\cG$ in $\cH$; see Fig. \ref{fg3-2}.
 Let $(x,\lam)\in \cK$. If $x=\theta$ then clearly  $(x,\lam)\in \mb{int}_\cH\,\cG$. (Here $\mb{int}_\cH \cV$ denotes the interior of $\cV\subset \cH$ in $\cH$.)  Thus we  assume $x\ne \theta$.  Then by \eqref{e3.12b} we have $$(x,\lam)\in \cF_\cH\subset \cF\subset \sP(\cC).$$ Thus by \eqref{e:3.6} we deduce that  $(x,\lam)\not\in\pa\cO$.  Since  $(x,\lam)\in \cF_\cH\subset \cO$, we also have
 $$
 (x,\lam)\not\in (\pa B_{\rho}\X\Lam)\backslash{\cO}.
 $$
Hence by  \eqref{e3.12g} one concludes  that $(x,\lam)\not\in \pa_\cH\cG$, which proves our claim.

Now the continuation property of the Conley index implies that $h(\Phi_\lam,\cK[\lam])$ remains constant on $\Lam$. In particular, we have
\be\label{e:3.10}h\(\Phi_{\a-\ve},\cK[\a-\ve]\)=h\(\Phi_{\b+\ve},\cK[\b+\ve]\).\ee

As is depicted in Fig. \ref{fg3-2}, by \eqref{e:3.8} and the choice of $\rho$ we have
\be\label{e:3.10a}\cO[\lam]\cap \ol B_{\rho}=\emp,\Hs \lam=\a-\ve,\,\,\b+\ve.\ee
{(In case $\cO[\lam]=\emptyset$ for $\lam=\a-\ve$ or $\b+\ve$, \eqref{e:3.10a} naturally holds.)}
Therefore by Lemma \ref{l:2.4.1}
\be\label{e3.21}
 h\(\Phi_{\lam},\cK[\lam]\)=h\(\Phi_{\lam},\theta\)\vee h\(\Phi_{\lam},\cF[\lam]\),\Hs \lam=\a-\ve,\,\,\b+\ve.
 \ee

\begin{figure}[h!]
  \centering
  \includegraphics[width=6.8cm]{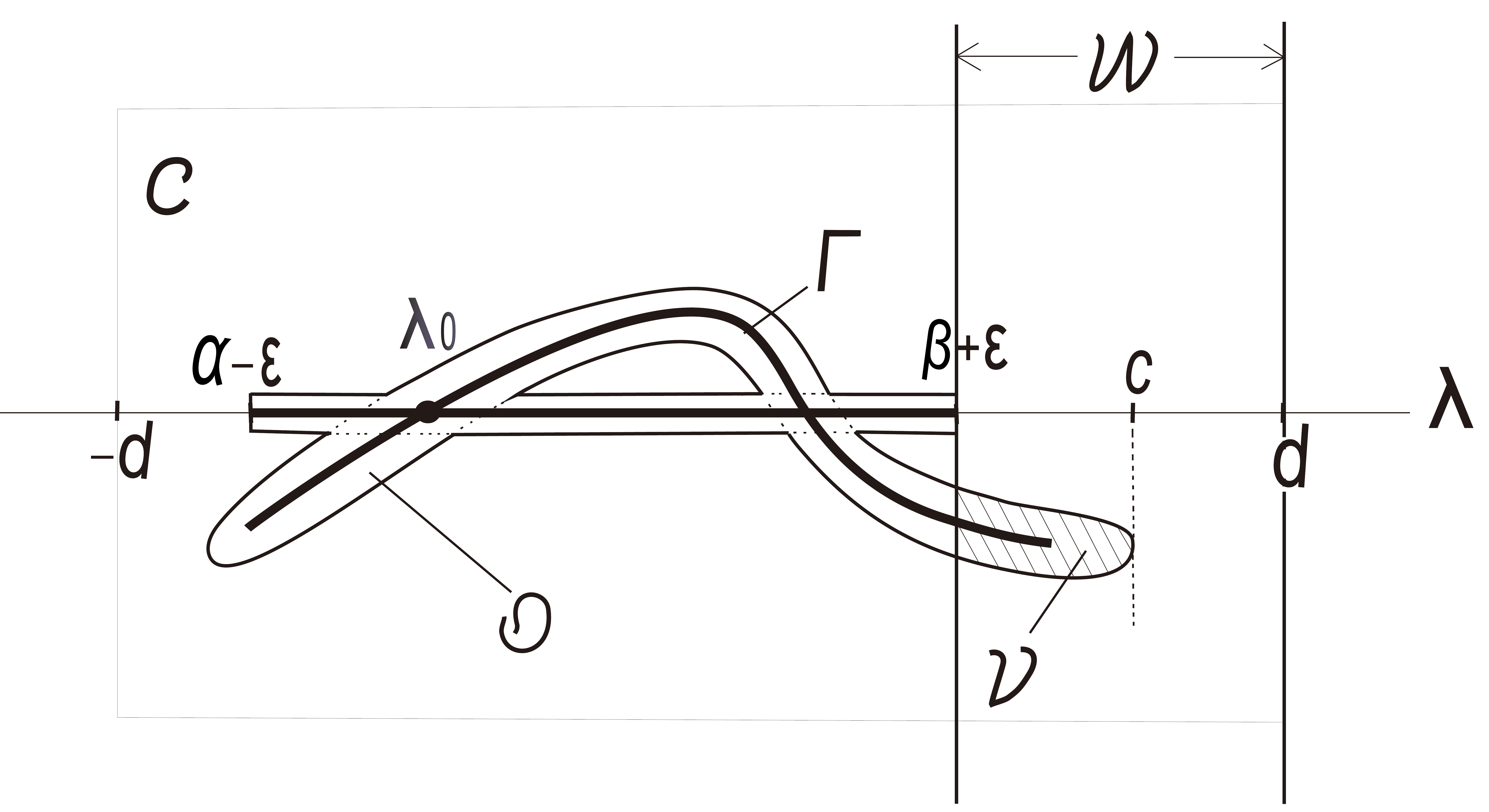}\\
  \caption{$\cV$  is an isolating neighborhood $\cF_\cW$ in  $\cW$}\label{fg3-1}
\end{figure}

 In what follows we check that
\be\label{e3.h} h\(\Phi_{\a-\ve},\cF[{\a-\ve}]\)=\ol0=h\(\Phi_{\b+\ve},\cF[{\b+\ve}]\).
 \ee 
{We only consider the case where $\lam=\b+\ve$. If $\cO[\b+\ve]=\emptyset$ then  \eqref{e3.h} clearly holds true. } So we assume $\cO[\b+\ve]\ne \emp$. Define
$$
c=\max\{\lam\geq \b+\ve|\,\,\,\cO[\lam]\ne\emp\}.
$$
Then $c<d$.
Let $$\cW=X\X[\b+\ve,d],\hs \cV=\cO\cap \cW;$$  see Fig. \ref{fg3-1}. One can easily see  that  $\pa_\cW \cV=\pa\cO\cap \cW$.
Since $\theta\notin \cV[\lam]$ for $\lam\in [\b+\ve,c]$, by a similar argument as in the verification of the isolating property of the domain $\cG$, we deduce that $\cV$ is an isolating neighborhood of $I(\~\Phi,\cV)=\cF\cap \cW:=\cF_\cW$ with respect to $\~\Phi$ restricted on $\cW$. Hence $\cF_\cW$ is an isolated compact invariant set of $\~\Phi$ in $\cW$.
By Theorem \ref{t:2.14} we have
\be\label{e3.j}h\(\Phi_{\lam},\cF_\cW[\lam]\)=h\(\Phi_{\lam},\cF[\lam]\)\equiv \mb{const}.,\Hs \lam\in[\b+\ve,\,d].\ee
Noting that  $\cF[\lam]=\emp$  for $\lam\in[c,d]$ (see Fig. \ref{fg3-1}), we  have  $h\(\Phi_{\lam},\cF[\lam]\)=\bar{0}$ for all $\lam\in [\b+\ve,\,d].$
In particular,
$h\(\Phi_{\b+\ve},\cF[{\b+\ve}]\)=\ol0.
$

Combining \eqref{e:3.10}, \eqref{e3.21} and \eqref{e3.h} together, we get
\be\label{e:3.11}h\(\Phi_{\a-\ve},\theta\)=h\(\Phi_{\b+\ve},\theta\).\ee
On the other hand, since $\lam_0$ is an essential bifurcation value, there is a number $\eta>0$ sufficiently small such that $
h\(\Phi_{\lam_0-\eta},\theta\)\ne h\(\Phi_{\lam_0+\eta},\theta\).
$ Because there are no other bifurcation values in the interval $[\a-\ve, \b+\ve]$ other than $\lam_0$, by the continuation property of the Conley index  we conclude that
$$
h\(\Phi_{\a-\ve},\theta\)=h\(\Phi_{\lam_0-\eta},\theta\)\ne h\(\Phi_{\lam_0+\eta},\theta\)=h\(\Phi_{\b+\ve},\theta\),
$$
which leads to a contradiction and completes the proof of Theorem \ref{t:3.2}. \,$\Box$

\subsection{A second global dynamic  bifurcation theorem}
\Vs

In Theorem \ref{t:3.2} we have shown that if the global bifurcation branch $\Gamma$ of an essential bifurcation point $(\theta,\lam_0)$ is bounded, then it necessarily meets the trivial  branch $\{\theta\}\X\R$ at another point $(\theta,\mu_0)$. However, it remains open whether  $(\theta,\mu_0)$ is a bifurcation point.  In this subsection we give an affirmative answer to this question  under some additional  reasonable assumptions.

We need  the following hypotheses:

 \benu\item[(H1)] Each bifurcation value $\lam$ is isolated.
\vs\item[(H2)]  If $\lam\not\in\Upsi$ then $h(\Phi_\lam,\theta)=\Sigma^p$ for some $p\geq 0$.
 \eenu

\br\label{r:4.2}Since $\Upsi$ is closed, under the hypothesis $(H1)$ one can easily verify  that each bounded interval contains only a finite number of bifurcation values $\lam\in\Upsi$.
\er

 \bt\label{t:3.3} $($A second global dynamic bifurcation theorem$)$ Assume  the hypotheses $(H1)-(H2)$.
 Let $(\theta,\lam_0)$ be an essential bifurcation point of $\{\Phi_\lam\}_{\lam\in\R}$.
 Then one of the following two alternatives occurs.
 \benu\item[$(1)$] The global bifurcation branch $\Gamma=\Gamma(\theta,\lam_0)$  is unbounded in $\sX$.
 \item[$(2)$] $\Gamma$ connects to  another bifurcation point  $(\theta,\lam_1)$ with
$\lam_1\ne\lam_0$.
\eenu
 \et

\noindent{\bf Proof.} We assume $\Gamma$ is bounded and prove that the assertion (2) holds.

\begin{figure}[h!]
  \centering
  \includegraphics[width=8.8cm]{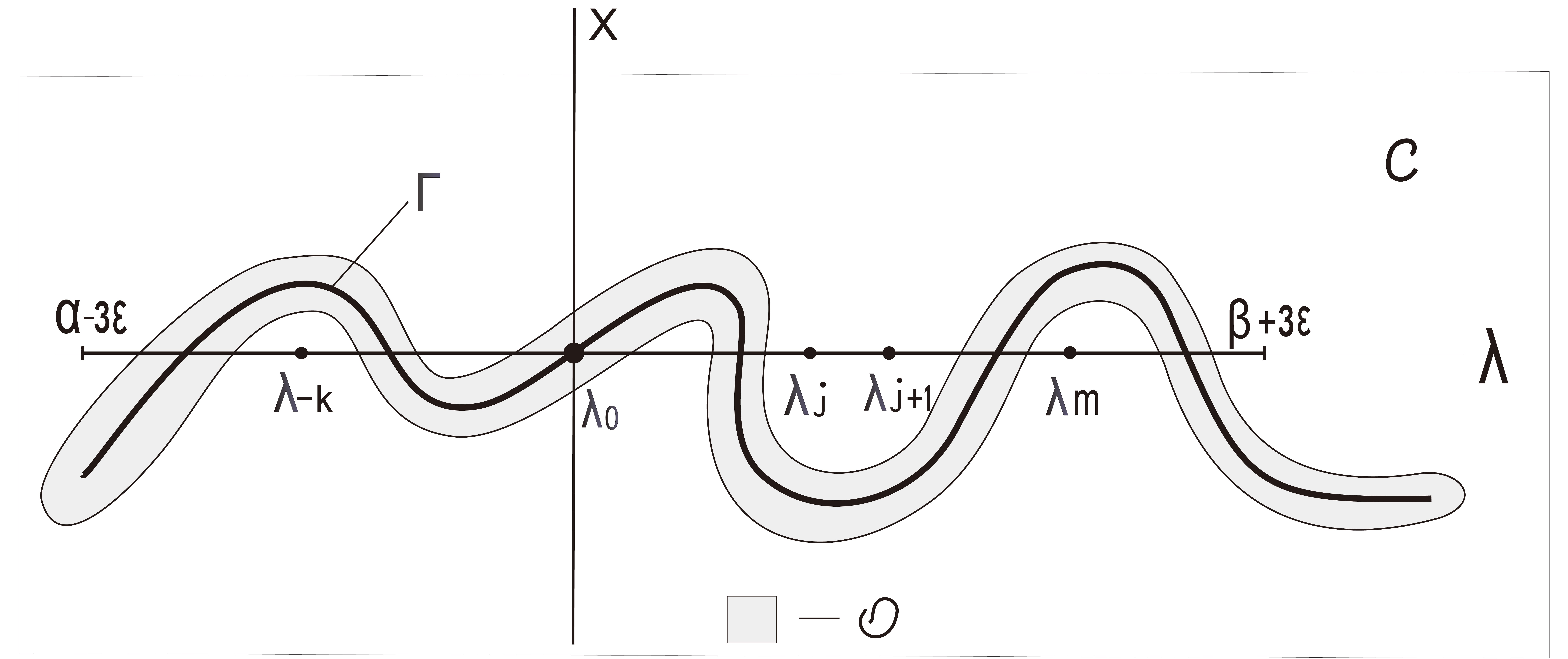}\\
  \caption{Distribution of ${\lam_j}'s$}
  \label{fg3-3}
\end{figure}

First, we infer from the proof of Theorem \ref{t:3.2} that $\Gamma$ is a connected compact invariant set of the skew-product flow $\~\Phi$. Let $J=\{\lam|\,\,(\theta,\lam)\in\Gamma\}$, and denote $$\a=\min\{\lam|\,\,\lam\in J\},\hs \b=\max\{\lam|\,\,\lam\in J\}.$$
Then $\a,\b\in J$. We infer from Theorem \ref{t:3.2} that the interval $[\a,\b]$ contains at least a bifurcation value $\lam_1\in\Upsi$ with $\lam_1\ne\lam_0$.
By Remark \ref{r:4.2} there are only a finite number of bifurcation values in $[\a,\b]$.
Set \be\label{e:U0}\ba{ll}[\a,\b]\cap \Upsi=\{\lam_j|\,\,-k\leq j\leq m\}:=\Upsi_0,\ea\ee where
$
\lam_{-k}<\cdots<\lam_{-1}<\lam_0<\lam_1<\cdots<\lam_m;
$
see Fig. \ref{fg3-3}. We show that there exists $\lam_j\ne\lam_0$ such that $\lam_j\in J$, thus proving what we desired.

We argue by contradiction and suppose that $\lam_j\not\in J$ for all $\lam_j\ne\lam_0$. Then
$$J\cap\Upsi=\{\lam_0\}.$$
By Remark \ref{r:4.2}  we can pick an $\ve>0$ such that
$$
[\a-4\ve,\,\b+4\ve]\cap\Upsi=[\a,\b]\cap\Upsi=\Upsi_0,\hs J_{4\ve}\cap\Upsi=\{\lam_0\},
$$
where $J_r$ denotes the $r$-neighborhood of $J$ (see \eqref{e:j3}). Therefore
\be\label{e:j4}
d(\lam_j,J)\geq 4\ve>0,\Hs -k\leq j\leq m,\,\,j\ne 0.
\ee

Taking two numbers $d,R>0$ so that $\ol\mB(\Gamma,1)\subset \ol B_R\X[-d,d]:=\cC$. Repeating the same argument leading to \eqref{e:3.8} and  \eqref{e:3.6} in the proof of Theorem \ref{t:3.2} with almost no modification, one can find a closed neighborhood $\cO$ of $\Gamma$ such that
\be\label{e:3.30}
d(\theta,\cO[\lam])\geq2\de_0>0,\Hs \A\,\lam\in {J_\ve}^{c}:=\R\sm J_\ve,
\ee
and
\be\label{e:3.29}
\mB(\pa\cO,2\de_\Gamma)\cap\sP(\cC)=\emp,
\ee
for some positive numbers $\de_0,\de_\Gamma>0$; see Fig. \ref{fg3-3}.

\begin{figure}[h!]
  \centering
  \includegraphics[width=8cm]{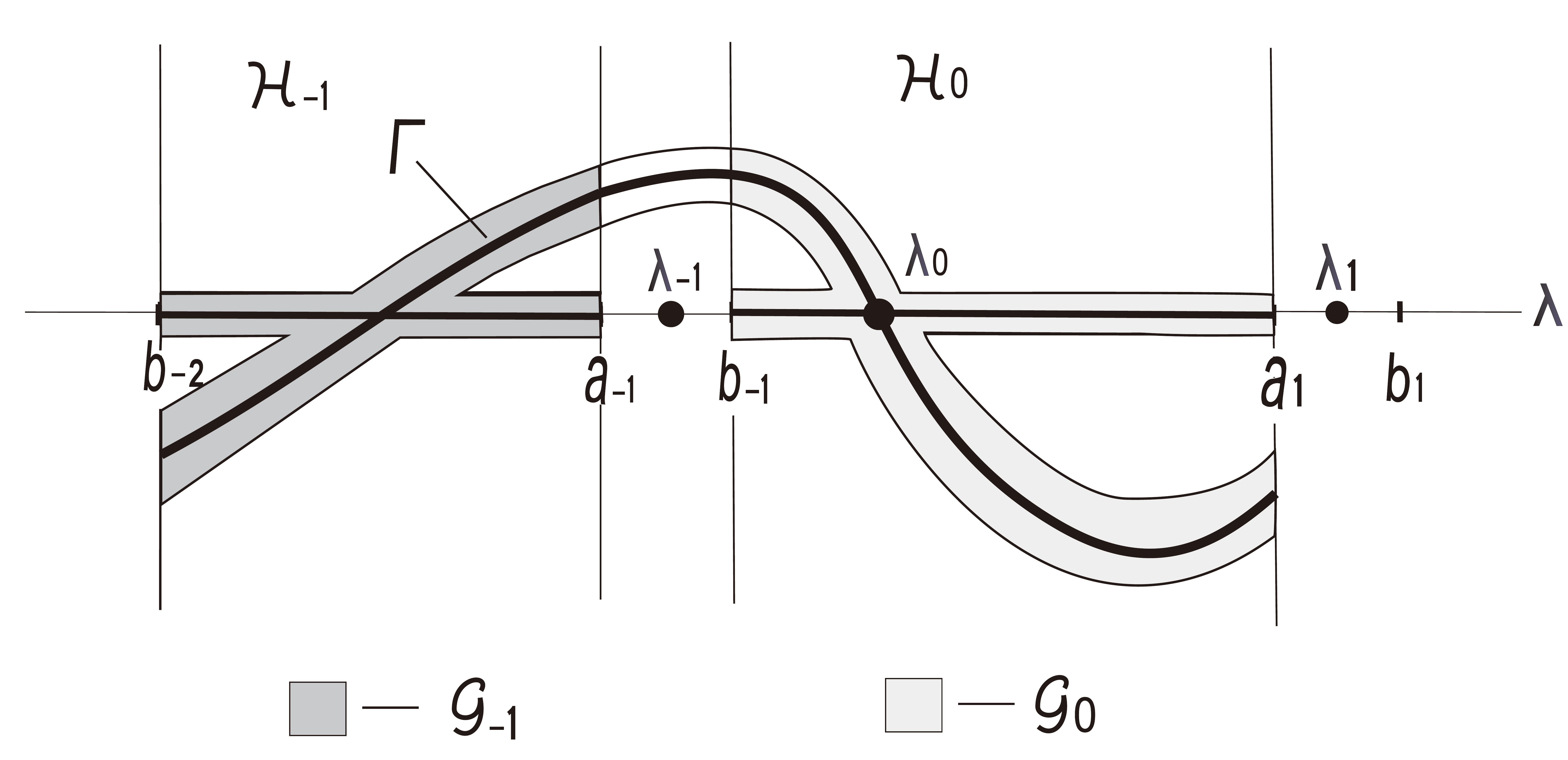}\\
  \caption{$(\theta,\lam_0)\in\cG_0$}
  \label{fg3-5}
\end{figure}

For notational simplicity, we assign $$b_{-k-1}=\a-3\ve,\hs a_{m+1}=\b+3\ve.$$ Write
$$
 \lam_j-\ve:=a_j,\hs \lam_j+\ve:=b_{j},\Hs -k\leq j\leq m,\,\,j\ne0.
$$
Note that $\lam_0\in[b_{-1},a_1]$; see Fig. \ref{fg3-5}. By \eqref{e:j4} it is clear that $a_i,b_{j}\notin \Upsilon$ for all  $-k\leq i\leq m+1$ and $-k-1\leq j\leq m$, $i,j\ne 0$; furthermore,
\be\label{e:3ab}
d(a_i,J)\geq 3\varepsilon,\hs d(b_{j},J)\geq 3\varepsilon.
\ee

\begin{figure}[h!]
  \centering
  \includegraphics[width=7.5cm]{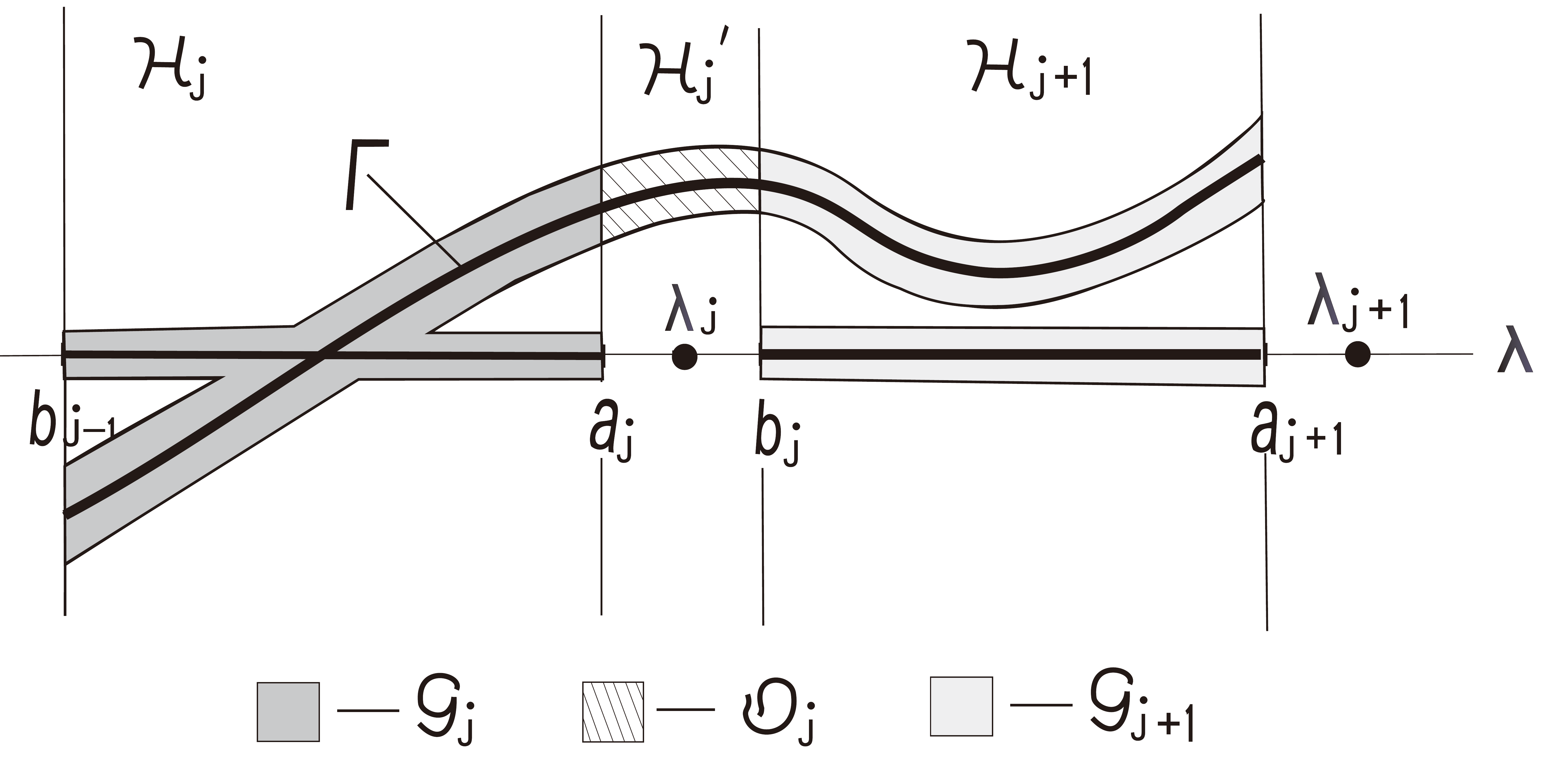}\\
  \caption{$\cO_j$ is isolating in $\cH_j'$}
  \label{fg3-4}
\end{figure}

{Let $\Lam_0=[b_{-1},\,a_1]$}, and let
$$
\Lam_j=[b_{j-1},\,a_{j}],\Hs  -k\leq j\leq m+1,\,\,\,j\ne0, 1.
$$
Set $\cH_j=X\X\Lam_j$. Following  the procedure  in the proof of Theorem \ref{t:3.2} below \eqref{e:3.8} with minor  modifications, one can choose a positive number $\rho<\de_0$ such that the following assertions hold:
\benu
\item[(I)]
 For each $-k\leq j\leq m+1$, $j\ne1$, the set $\ba{ll}\cG_j=\(\cO\cap \cH_j\)\Cup \(\ol B_\rho\X \Lam_j\)\ea$ is an isolating neighborhood of $\~\Phi$ restricted on $\cH_j$ with
          \be\label{e:I1}\ba{ll} I(\~\Phi,\cG_j)=\(\cF\cap \cH_j\)\Cup (\{\theta\}\X\Lam_j),\ea\ee
         {where $\cF=I(\~\Phi,\cO)$};  see Figures  \ref{fg3-5} and \ref{fg3-4}.
\item[(II)] For each  $-k\leq j\leq m$, $j\neq 0$, the set $\cO_j=\cO\cap \cH_j'$ is an isolating neighborhood of $\~\Phi$ restricted on $\cH_j'=X\X [a_j,b_{j}]$ with $I(\~\Phi,\cO_j)=\cF\cap \cO_j$;  see Fig. \ref{fg3-4}.
\eenu
\vs

By \eqref{e:3.30}, \eqref{e:3ab} and the choice of $\rho$ it is easy to see that  $$\cO[a_i]\cap \ol B_\rho=\emp=\cO[b_j]\cap \ol B_\rho$$  for all $a_i$ and $b_j$. Thus by the assertion (I), Lemmas \ref{l:2.4.1}  and Theorem \ref{t:2.14} we deduce that
\begin{equation}\label{e:3.31b}
      h(\Phi_{b_{-1}},\cF[b_{-1}])\vee h(\Phi_{b_{-1}},\theta)=h(\Phi_{a_1},\cF[a_1])\vee h(\Phi_{a_1},\theta),
\end{equation}
and
\begin{equation}\label{e:3.31}
     h(\Phi_{b_{j-1}},\cF[b_{j-1}])\vee h(\Phi_{b_{j-1}},\theta)= h(\Phi_{a_j},\cF[a_j])\vee h(\Phi_{a_j},\theta)
\end{equation}
for all $-k\leq j\leq m+1$, $j\ne 0,1$. By the assertion (II) and Theorem \ref{t:2.14} we have
\be\label{e:3.32}
h(\Phi_{a_j},\cF[a_j])=h(\Phi_{b_{j}},\cF[b_{j}])
\ee
for $-k\leq j\leq m$, $j\ne 0$.

In what follows we first check that
\begin{equation}\label{e:3.37}
\begin{split}
&h(\Phi_{a_1},\cF[a_1])=\ol0.
\end{split}
 \end{equation}
Recall that $a_{m+1}=\b+3\ve$. Using  a similar argument as in the verification of  \eqref{e3.h} it can be shown that
 \begin{equation}\label{e:3.34}
 h(\Phi_{a_{m+1}},\cF[a_{m+1}])=\ol0.
 \end{equation}
If $m=0$ (i.e., the interval $(\lam_0,b+3\ve]$ contains no  bifurcation values) then  we are done. Thus we assume that $m\geq 1$.

As $b_{m},a_{m+1}\notin\Upsilon$, by the hypothesis $\text{(H2)}$
we find that
\begin{equation*}
h(\Phi_{b_{m}},\theta)=\Sigma^{p},\hs h(\Phi_{a_{m+1}},\theta)=\Sigma^{q}
\end{equation*}
for some nonnegative integers $p$ and $q$. On the other hand, by \eqref{e:3.31} we have
\begin{equation*}\ba{ll}
      h(\Phi_{b_{m}},\cF[b_{m}])\vee h(\Phi_{b_{m}},\theta)&=h(\Phi_{a_{m+1}},\cF[a_{m+1}])\vee h(\Phi_{a_{m+1}},\theta)\\[1ex]
      &=(\mb{by \eqref{e:3.34}})\\[1ex]
      &=\ol0\vee h(\Phi_{a_{m+1}},\theta)=\Sigma^q.\ea
\end{equation*}
Therefore by  Lemma \ref{l£»2.2-1} we necessarily have
 $ h(\Phi_{b_{m}},\cF[b_{m}])=\ol0.
$
Further by \eqref{e:3.32} one concludes  that
\begin{equation}\label{e:3.36}
 h(\Phi_{a_m},\cF[a_m])=\ol0.
 \end{equation}

If $m=1$ then we are done. Otherwise one can repeat the above argument with $a_{m}$, $b_{m-1}$ and $a_{m-1}$ in place of $a_{m+1},\,b_m$ and $a_{m}$ respectively to find that
 $$ h(\Phi_{b_{m-1}},\cF[b_{m-1}])=\ol0,\hs h(\Phi_{a_{m-1}},\cF[a_{m-1}])=\ol0.$$ Proceeding this procedure   we finally conclude the validity of \eqref{e:3.37}.

A parallel argument as above applies to show that
\begin{equation}\label{e:3.38}
\begin{split}
&h(\Phi_{b_{-1}},\cF[b_{-1}])=\ol0.
\end{split}
 \end{equation}

Combining  \eqref{e:3.31b},  \eqref{e:3.37} and  \eqref{e:3.38} it yields
$$
h(\Phi_{b_{-1}},\theta)=h(\Phi_{a_1},\theta). $$
However,  $\lam_0$ is an essential bifurcation value; moreover, it is the unique bifurcation value in the interval $[b_{-1}, a_1]$.   By the definition of an essential bifurcation value and the continuation property of the Conley index, it is easy to deduce that
$
h(\Phi_{b_{-1}},\theta)\ne h(\Phi_{a_1},\theta),
$
which leads to a contradiction and completes the proof of Theorem \ref{t:3.3}. \,$\Box$

\subsection{A third global dynamic  bifurcation theorem}
\Vs

Finally, let us give a third  global dynamic bifurcation theorem in which there is only one possibility, that is, the bifurcation branch is unbounded.
For this purpose, we need to impose on the Conley index along the trivial equilibrium  branch $\{\theta\}\X \R$ a stronger  condition:

 \benu\item[(H3)] For any compact interval $[a,b]\subset\R$ with $a,b\notin \Upsi$, we have
 $$
 h(\Phi_a,\theta)\ne h(\Phi_b,\theta)
 $$
 whenever  $(a,b)\cap\Upsi\ne\emp$.
 \eenu
 At first glance, this requirement seems to be quite   restrictive. However,  due to the {\em nonnegativity} of the Conley index, it is naturally fulfilled by a large number of examples from applications (as we will see in Section \ref{S5}).
 \bt\label{t:3.7} $($A third global dynamic bifurcation theorem$)$  Assume  the hypotheses $($H1$)-($H3$)$.  Then for any  bifurcation point  $(\theta,\lam_0)$, the global dynamic bifurcation branch $\Gamma=\Gamma(\theta,\lam_0)$  is unbounded in $\sX$
 \et

\noindent{\bf Proof.}
We argue by contradiction and suppose that $\Gamma$ is bounded. Then as in the proof of Theorem \ref{t:3.2}, we deduce that $\Gamma$ is a connected compact invariant set of $\~\Phi$. Set $J=\{\lam|\,\,(\theta,\lam)\in\Gamma\}$, and denote $$\a=\min\{\lam|\,\,\lam\in J\},\hs \b=\max\{\lam|\,\,\lam\in J\}.$$
Clearly $\a,\b\in J$.
 Let $$\Upsi_0=\{\lam\in \Upsi|\,\,\lam\in [\a,\b]\sm J\}.$$
(One should distinguish the set $\Upsi_0$ defined above  from the one  in \eqref{e:U0}.)
By Remark \ref{r:4.2} we may write
$ \Upsi_0=\{\lam_j|\,\,-k\leq j\leq m,\,\,j\ne 0\},$ where
$$
\lam_{-k}<\cdots<\lam_{-1}<\lam_0<\lam_1<\cdots<\lam_m;
$$
see Fig. \ref{fg3-3} for the distribution of ${\lam_j}'s$.

The remaining part of the  argument  is almost the same as in the proof of Theorem \ref{t:3.3} except that, instead of assuming $\lam_0$ is an essential bifurcation value, we employ a stronger assumption (H3).

By (H1) we can pick an $\ve>0$ such that
$$
[\a-4\ve,\,\b+4\ve]\cap\Upsi=[\a,\b]\cap\Upsi,\hs J_{4\ve}\cap\Upsi_0=\emp.
$$
Then
$$
d(\lam_j,J)\geq 4\ve>0,\Hs \lam_j\in\Upsi_0.
$$
Take two numbers $d,R>0$ so that $\ol\mB(\Gamma,1)\subset \cC:=\ol B_R\X[-d,d]$. Then as in the proof of Theorem \ref{t:3.3},  one can find a closed neighborhood $\cO$ of $\Gamma$ such that \eqref{e:3.30} and \eqref{e:3.29} remain valid for some  $\de_0,\de_\Gamma>0$.

For notational simplicity, we assign $$b_{-k-1}=\a-3\ve,\hs a_{m+1}=\b+3\ve,$$ and write
$$
 \lam_j-\ve:=a_j,\hs \lam_j+\ve:=b_{j},\Hs -k\leq j\leq m,\,\,j\ne0.
$$
Clearly $$
d(a_i,J)\geq 3\varepsilon>0,~~\text{and}~~d(b_{j},J)\geq 3\varepsilon>0,
$$
  for all  $-k\leq i\leq m+1$ and $-k-1\leq j\leq m$, $i,j\ne 0$.

{Let $\Lam_0=[b_{-1},\,a_1]$}, and let
$$
\Lam_j=[b_{j-1},\,a_{j}],\Hs  -k\leq j\leq m+1,\,\,\,j\ne0, 1.
$$
Following  the procedure  in the proof of Theorem \ref{t:3.2} below \eqref{e:3.8} with minor  modifications, one can choose a positive number $\rho<\de_0$ such that the assertions (I) and (II) in the proof of Theorem \ref{t:3.3} hold true.
Therefore by the continuation property of the Conley index we see that \eqref{e:3.31b}, \eqref{e:3.31} and \eqref{e:3.32} remain valid.
Now repeating the same argument as in the proof of Theorem \ref{t:3.3}, one finds that
$$h(\Phi_{a_1},\cF[a_1])=\ol0=h(\Phi_{b_{-1}},\cF[b_{-1}]).
$$
Further by \eqref{e:3.31b} we deduce that
$
h(\Phi_{b_{-1}},\theta)=h(\Phi_{a_1},\theta). $
On the other hand, since $\lam_0\in [b_{-1},\,a_1]$ and $b_{-1},a_1\notin\Upsi$, by (H3) we have
$
h(\Phi_{b_{-1}},\theta)\ne h(\Phi_{a_1},\theta).
$
This leads to a contradiction. \,$\Box$

\section{Global  Bifurcation of Evolution Equations}\label{S5}
Our first example of applications of the abstract bifurcation theorems given in Section 4 concerns  the evolution  equation
\be\label{e:3.1}
u_t+A u=f_\lam(u)
\ee in a Banach space $X$,
where  $A$ is a sectorial operator in $X$ with compact resolvent,  and $f_\lam(u)$ is a locally Lipschitz continuous  mapping from $X^\alpha\X \R$ to $X$ for some $0\leq\alpha<1$. (We refer the reader to Henry \cite{Henry} for a basic  theory on sectorial operators and fractional powers of Banach spaces.)

It is well known that under the above hypotheses, the initial value problem of \eqref{e:3.1} is well-posed in $X^\a$. That is, for each  $u_0\in X^\a$ the equation \eqref{e:3.1}  has a unique strong solution $u(t)$ in $X^\a$ with $u(0)=u_0$ on a maximal existence interval $[0,T_{u_0})$ $($see e.g. \cite[Theorem 3.3.3]{Henry}$)$.
Denote $\Phi_\lam$ the local semiflow generated by the initial value problem of \eqref{e:3.1} on $X^\a$, namely, given $u_0\in X^\a$, $u(t)=\Phi_\lam(t)u_0$ is precisely the solution of \eqref{e:3.1} with initial value $u(0)=u_0$.
\vs
Set $\sX=X^\a\X\R$, and let  $\~\Phi$ be the {\em skew-product flow} of the family $\{\Phi_\lam\}_{\lam\in \R}$ on $\sX$.
Then by standard argument (see e.g.\,\,\cite[Chap. 3, Theorem 3.3.6]{Henry} or \cite[Chap. I, Theorem 4.4]{Ryba}), one can easily verify  that $\~\Phi$ is {\em asymptotically compact}, i.e., $\~\Phi$ satisfies the hypothesis  (AC) in Section 2.

Suppose $f_\lam(0)=0$ for all $\lam\in\R$, hence $u=0$ is always a trivial equilibrium solution of \eqref{e:3.1}.
We also assume that $f_\lam(u)$ is differentiable  in $u$ with $\partial_uf_\lam(u)$ being continuous in $(u,\lam)$.
Let $$L_\lam=A-\partial_u f_\lam(0).$$ Then  $L_\lam$  is a sectorial operator in $X$ with {\em compact resolvent}; See Proposition \ref{p:7.1} in the Appendix.

 \eqref{e:3.1} can be rewritten as
\be\label{e3.2}
u_t+L_\lam u=g_\lam(u),
\ee
where $g_\lam(u)=f_\lam(u)-\partial_u f_\lam(0)\,u$.
Denote
$$\~\Upsilon=\{\lam|\,\,\mb{Re}\,\mu=0\mb{ for some }\mu\in\sig(L_\lam)\},$$
where $\sig(L_\lam)$ is the spectrum of $L_\lam$. One can easily  see that $\~\Upsi$ is closed in $\R$. If $\lam \notin \~\Upsilon$ then the trivial equilibrium solution $\theta=0$  is hyperbolic, and no bifurcation occurs near the  point $(0,\lam)$. Furthermore, we infer from Rybakowski \cite[Chapter II, Theorem 3.5]{Ryba} that  $$h(\Phi_\lam,0)=\Sigma^p$$ for some $p\geq0$.

Let $\Upsi$ be the set of bifurcation values $\lam\in\R$ of the system along the trivial equilibrium solution branch $\{0\}\X\R$. Then $\Upsi\subset \~\Upsi$.
\bl\label{l:5.1} Let $\lam_*\in\~\Upsi\sm\Upsi$. If $\lam_*$ is isolated  in $\~\Upsi$, then $h(\Phi_{\lam_*},0)=\Sigma^p$ for some $p\geq0$.
\el

\noindent{\bf Proof.} Since $\lam_*$ is isolated  in $\~\Upsi$, there exists $\ve>0$ such that $[\lam_*-\ve,\lam_*+\ve]\cap \~\Upsi=\{\lam_*\}$. Consequently
$[\lam_*-\ve,\lam_*+\ve]\cap \Upsi=\emp$. Hence one can easily deduce that there is a closed neighborhood  $N$ of $0$ in $X^\a$ such that  $N$ is an isolating neighborhood of the invariant set $S_0=\{0\}$ of $\Phi_\lam$ for all $\lam\in[\lam_*-\ve,\lam_*+\ve]$.
The continuation property of the Conley index then implies that
\be\label{e:l5.1}
h(\Phi_\lam,0)\equiv \mb{const.}, \Hs \lam\in [\lam_*-\ve,\lam_*+\ve].
\ee
 Because $\lam_1=\lam_*-\ve\notin \~\Upsi$, we have
$h(\Phi_{\lam_1},0)=\Sigma^p$ for some $p\geq0$. Thus by \eqref{e:l5.1} we conclude that $h(\Phi_{\lam_*},0)=h(\Phi_{\lam_1},0)=\Sigma^p$. \,$\Box$
\Vs

Assume  $\lam_0\in \~\Upsi$ satisfies the following hypothesis:
\benu
\item[(A1)] There exists $\ve>0$ such that for $\lam\in[\lam_0-\ve,\,\lam_0+\ve]$, the spectrum  $\sig(L_\lam)$ has a decomposition $\ba{ll}\sig(L_\lam)=\sig_\lam^1\cup\sig_\lam^2\cup\sig_\lam^3\ea $ with $\sig_\lam^2\ne\emp$ such  that
\begin{equation}\label{e:A1}
\sig_\lam^1\subset \mathbb{C}^-,\hs \sig_\lam^3\subset \mathbb{C^+},
\end{equation}
 where $\mathbb{C}^\pm=\{z\in \mathbb{C}|\,\,\pm\mb{Re}\,z>0\}$, and
\begin{equation}\label{e:H2}
\sig_\lam^2\subset \mathbb{C^+}\,\,(\mb{if }\,\lam<\lam_0),\hs \sig_\lam^2\subset \mathbb{C}^-\,\,(\mb{if }\,\lam>\lam_0).
\end{equation}
\eenu
Since $L_\lam$ has compact resolvent,   $\sig_\lam^1$ and $\sig_\lam^2$ necessarily consist of a finite number of eigenvalues of $L_\lam$. As $\sig_\lam^2\ne\emp$, we deduce  by \eqref{e:A1}, \eqref{e:H2} and \cite[Chapter II, Theorem 3.1]{Ryba}  that the Conley index $h(\Phi_\lam,0)$ changes as $\lam$ crosses $\lam_0$ (see also \cite[Section 4]{LW}). Hence  $\lam_0$ is an essential bifurcation value.

\vs

By virtue of Lemma \ref{l:5.1},  if each point in $\~\Upsi$ is isolated then the hypothesis (H2) in  Theorem \ref{t:3.3} is satisfied. As a direct application of Theorem \ref{t:3.3}, we immediately obtain the following global dynamic bifurcation result.

 \bt\label{gbt1} Suppose that each point in $\~\Upsi$ is isolated. Assume that $\lam_0\in\~\Upsi$ satisfies ${(A1)}$.
 Let $\Gamma=\Gamma(0,\lam_0)$ be the global dynamic bifurcation branch of $(0,\lam_0)$. Then either  $\Gamma$  is unbounded, or $\Gamma$ meets  another bifurcation point  $(0,\lam_1)$ with $\lam_1\ne \lam_0$.
 \et

 A particular but important case is the system
 \be\label{e:3.1c}
u_t+A u=\lam u+f(u).
\ee
     Since  $A-f'(0)$ has compact resolvent, $\sig(A-f'(0))$ consists of eigenvalues $\mu_k$ ($k=1,2,\cdots$) with $\mb{Re}\,\mu_k\ra +\8$. Thus  one easily sees  that
 $$
 \Upsi=\~\Upsi=\{\lam_k:=\mb{Re}\,\mu_k|\,\,\,k\geq 1\}.
 $$

Let $-\8<a<b<\8$.  Assume $a,b\notin \Upsi$  and that $(a,b)\cap\Upsi\ne\emp$.  Then we infer from \cite[Chap. II, Theorem 3.5]{Ryba} that
 $
 h(\Phi_a,0)=\Sigma^p$ for some $p\geq 0$, whereas  $h(\Phi_b,0)=\Sigma^{p+r}, $
 where $r>0$ is the sum of the algebraic multiplicities of the eigenvalues $\mu_k$ with $\mb{Re}\,\mu_k\in (a,b)$. Hence we see that the hypothesis (H3) in Theorem \ref{t:3.7} is also fulfilled. As a result, we have

  \bt\label{gbt2}  For each $k\geq 1$, the global dynamic bifurcation branch $\Gamma$ of $(0,\lam_k)$ $($with respect to the system \eqref{e:3.1c}$)$ is unbounded.
 \et

Theorem \ref{gbt2} seems to be somewhat different from the known ones  given in the literature, and may help us have a deeper understanding of the dynamics of nonlinear evolution equations.

\section{Applications to Elliptic Problems}

 Let $\Omega\subset \R^n$ ($n\geq 3$) be a bounded open domain.
 Consider the  equation:
\be\label{e:4.1}
-\Delta u=f_\lam(u),\Hs x\in\Omega
\ee
associated with the homogeneous Dirichlet boundary condition, where $\lam\in\R$, and $f_\lam(s)$
  is a continuous function which is also assumed  to be differentiable in $s$ with $f'_\lam(s)$ being continuous in $(s,\lam)$.

We always assume $f_\lam(0)=0$ for $\lam\in\R$, so that $\{0\}\X\R$ is a trivial solution branch of \eqref{e:4.1}.

\subsection{Mathematical setting and the main result}
\Vs

Let $H=L^2(\Omega)$, and $V=H^1_0(\Omega)$. Denote by  $|\.|$ the usual norm of $L^2(\Omega)$, and define the norm $||\.||$ of $H^1_0(\Omega)$ as follows:
$$||u||=\(\int_\Omega |\nabla u|^2 dx\)^\frac{1}{2}.$$
We also use $|\.|_q$ to denote the norm of $L^q(\Omega)\ (q\geq 1)$.


We will impose on $f_\lam$ the following conditions:

\benu
\item[$(\mb{f}_1)$] There exists $1\leq p<\frac{n+2}{n-2}$ such that for any bounded interval $\Lam$,
$$|f_\lam'(s)|\leq a_1+a_2|s|^{p-1},\Hs \forall s\in \R,\,\,\lam\in\Lam$$
 for some positive constants $a_1$ and $a_2$.

\item[$(\mb{f}_2)$]There exists $\mu >2$ such that for any bounded interval $\Lam$ and $\ve>0$,
$$sf_\lam(s)\geq \mu F_\lam(s)-\varepsilon|s|^2-C_\varepsilon,\Hs \forall s\in \R,\,\,\lam\in\Lam$$
for some $C_\ve>0$, where $F_\lam(s)=\int_0^s f_\lam(s)ds$.
\end{enumerate}
Note that  $(\mb{f}_1)$ implies   that for any bounded interval $\Lam$,
\benu
\item[$(\widetilde{{\mb{f}}_1})$] there exist positive constants $a_3$ and $a_4$ such that
$$|f_\lam(s)|\leq a_3+a_4|s|^p,\Hs \forall s\in \R,\,\,\lam\in\Lam.$$
\eenu

Denote  $\mu_k$ $({k\geq 1})$ the distinct  eigenvalues of $A=-\Delta$ subjects to the homogeneous Dirichlet boundary condition,
$
0<\mu_1<\mu_2<\cdots<\mu_k<\cdots.
$

Our main results are summarized in the following theorem.

\bt\label{t:4.1} Assume $f$ satisfies $\(\mb{f}_1\)-\(\mb{f}_2\)$, and that  $\b(\lam)=f_\lam'(s)|_{s=0}$ is strictly increasing in $\lam$. Set
\be\label{e:Up}
 \Upsilon=\{\gam_k\in \R|\,\,\b(\gam_k)=\mu_k,\,\,k\geq 1\}.
 \ee
  Then for each $\gam_k$, there is an interval $\Lambda$ with $\gam_k\in\Lambda$  such that  $(\ref{e:4.1})$ has at least a nontrivial solution $u_\lam\ne 0$ for all $\lam \in \Lambda\setminus\Upsilon$. Moreover,  one of the  following two assertions holds:

\begin{enumerate}
\item[$(1)$] There is a bounded  sequence $\lam_m\in \Lambda$ such that  $||u_{\lam_m}||\rightarrow \infty$ as $m\rightarrow\8$.
\vs
\item[$(2)$] $\Lambda$ contains either the interval  $(-\8,\gam_k]$ or the interval $[\gam_k,\8)$.

\end{enumerate}
\et
Before proving  Theorem \ref{t:4.1}, we first give two illustration examples.
\Vs
\noindent{\bf Example 6.1.}
Consider the  equation:
\be\label{e:4.28}
-\Delta u=\lam g(u)+f(u),\Hs \mb{ in } \,\Omega,
\ee
where $g$ and $f$ are $C^1$ functions with $g(0)=0=f(0)$ and $g'(0)>0$. (We need not assume $f(s)=o(|s|)$ as $|s|\ra 0$.) Suppose $f$ satisfies $\(\mb{f}_1\)-\(\mb{f}_2\)$ with $f_\lam$ therein replaced by $f$, and that
$g$ is sublinear, i.e., there exist $\sigma\in[0,1)$ and $c_1$, $c_2>0$ such that
\be\label{e:g1}|g(s)|\leq c_1+c_2|s|^\sigma,\Hs s\in\R.\ee
Then $f_\lam(s)=\lam g(s)+f(s)$ satisfies $\(\mb{f}_1\)-\(\mb{f}_2\)$.

Indeed, by \eqref{e:g1}  one trivially verifies that  $f_\lam$  satisfies $\(\mb{f}_1\)$.

Write $g_\lam(s)=\lam g(s)$, and set $G_\lam(s)=\int_0^s g_\lam(t)dt$.  Let $\Lambda_0=[-a,a]$ $(a>0$), and let $\ve>0$ be given arbitrary.  Then by \eqref{e:g1} we have
\be\label{e:g2}|G_\lam(s)|\leq c_3|s|+c_4|s|^{\sigma+1},\Hs s\in\R,~~\lam\in\Lambda_0.\ee
Let  $\mu>2$ be the number given in (f$_2$). By \eqref{e:g1} and \eqref{e:g2} it is trivial to see that
$$
\left({\mu G_\lam(s)-sg_\lam(s)}\right)/{s^2}\rightarrow 0\hs \text{as}~|s|\rightarrow\8
$$  uniformly with respect to $\lam\in\Lam_0$.
Thus  there exists $M_\varepsilon>0$ such that
\begin{equation}\label{e:4.29}
-\varepsilon s^2\leq \mu G_\lam(s)-sg_\lam(s) \leq\varepsilon s^2,\Hs |s|>M_\varepsilon,\,\,\lam\in\Lam_0.
\end{equation}
Take a number $C_\ve>0$ such that
$
|\mu G_\lam(s)-sg_\lam(s)|\leq C_\varepsilon$ for all $|s|\leq M_\varepsilon$.
Then by \eqref{e:4.29} we have
\begin{equation}\label{e:4.31}
\mu G_\lam(s)-sg_\lam(s) \leq\varepsilon s^2+C_\varepsilon,~~\text{for all}~s\in\R.
\end{equation}
Hence we see that  $g_\lam$ satisfies $\(\mb{f}_2\)$ with $f_\lam$ therein replaced by $g_\lam$. Combining this with the assumption on $f$, one immediately concludes that $f_\lam(s)=\lam g(s)+f(s)$ satisfies $\(\mb{f}_2\)$.

\Vs\noindent {\bf Example 6.2.}
Consider the  equation:
\be\label{e:4.32}
-\Delta u=\lam u+\a|u|^{p-1}u+ \beta |u|^{q-1}u,\Hs \mb{ in } \,\Omega,
\ee
where $1<q<p<(n+2)/(n-2)$, and $\a,\b\in\R$ are constants with $\a\ne0$.

To apply Theorem $\ref{t:4.1}$, we only need to  check that the function $f_\lam(s)=\lam s+\a|s|^{p-1}s+\beta |s|^{q-1}s$ satisfies $\(\mb{f}_2\)$.
Set
$$F_\lam(s)=\int_0^s f_\lam(t)dt=\frac{\lam}{2}s^2+\frac{\a}{p+1}|s|^{p+1}+\frac{\beta}{q+1}|s|^{q+1},$$
and let  $\Lambda_0=[-a,a]$ $(a>0)$.
For $\lam\in\Lambda_0$, by simple  calculation we obtain that
\begin{equation*}
\begin{split}
\frac{\mu F_\lam(s)-sf_\lam(s)}{s^2}=&\lam\(\frac{\mu}{2}-1\)+\a\(\frac{\mu}{p+1}-1\)|s|^{p-1}\\&+\beta\(\frac{\mu}{q+1}-1\)|s|^{q-1}.
\end{split}
\end{equation*}
Since $p-1>q-1>0$, it is easy to see that
\vs
\begin{itemize}
  \item[(i)] if $\a>0$, then for $2<\mu<p+1$ we have
  $
 \({\mu F_\lam(s)-sf_\lam(s)}\)/{s^2}\rightarrow -\8
 $
  uniformly with respect to $\lam\in\Lam_0$, and
  \vs
  \item[(ii)] if $\a<0$, then for $\mu>p+1$ we have
  $
 \({\mu F_\lam(s)-sf_\lam(s)}\)/{s^2}\rightarrow -\8
 $
   uniformly with respect to $\lam\in\Lam_0$.
\end{itemize}
\vs
Thus in any case, one can pick two positive numbers $\mu>2$ and $R>0$ such that
\begin{equation*}
 \({\mu F_\lam(s)-sf_\lam(s)}\)/{s^2}<0,\Hs |s|>R, \,\,\lam\in\Lam_0.
  \end{equation*}
Hence
\begin{equation}\label{e:4.35}
\mu F_\lam(s)<sf_\lam(s),\Hs |s|>R, \,\,\lam\in\Lam_0,
\end{equation}
from which it can be easily seen  that  $f_\lam$ satisfies $\(\mb{f}_2\)$.
\vs

By virtue of Theorem \ref{t:4.1} one immediately obtains some global features on the bifurcation of the equations \eqref{e:4.28} and  \eqref{e:4.32}. For instance,  for the equation  \eqref{e:4.32}  the  fundamental  results summarized in Proposition \ref{t:1.1} hold true.

\subsection{Nonclassical parabolic flow}
\Vs

The basic idea  to prove Theorem \ref{t:4.1} is to regard  \eqref{e:4.1} as the stationary problem of  the nonclassical parabolic problem:
\be\label{e:4.3}\left\{\ba{lll}
u_t-\Delta u_t-\Delta u=f_\lam(u),\Hs &(x,t)\in \Omega\times \R^+,\\[1ex]
u(x,t)=0,\Hs &(x,t)\in \partial\Omega\times \R^+\ea\right.
\ee
and  apply global dynamic bifurcation theorems. Nonclassical parabolic equations have rich physical background and have attracted much interest in recent years; see e.g. \cite{Afan, SY,WLZ} and references cited therein.

Let us think of $A=-\De$ as an operator from $V$ to $V^*$:
 $$<Au,v>=\int_\Omega\nabla u\cdot\nabla vdx,~~~~~\forall~u, v\in V,$$
 where $V^*=H^{-1}_0(\Omega)$ and $<\cdot,\cdot>$ is the dual between $V$ and $V^*$.   

 We  use the same notation $f_\lam(\.)$ to denote the {\em Nemitski operator} given by the nonlinearity  $f_\lam(s)$.  Then for each $u\in V$, by $(\tilde{\mb{f}}_1)$  we have $$f_\lam(u)\in L^{(p+1)/p}(\Omega)=(L^{p+1}(\Omega))^*\subset V^*,$$
 and
 $$
 <f_\lam(u),v>=\int_\Omega f_\lam(u)vdx,\Hs v\in V.
 $$
 (\ref{e:4.3}) can be transformed into an abstract equation in $V$:
 \begin{equation}\label{AE0}
u_t+Au_t+Au=f_\lam(u),
\end{equation}
 or equivalently
\begin{equation}\label{AE}
u_t+Lu=g_\lam(u),
\end{equation}
where
$$L=(1+A)^{-1}A,~~~\text{and}~~g_\lam=(1+A)^{-1}f_\lam.$$
It is easy to deduce that $(1+A)^{-1}:V^*\rightarrow V$ is compact and that $L:V\rightarrow V$ is a bounded linear operator.
Hence $g_\lam$ is a nonlinear operator from $V$ to $V$. Therefore \eqref{AE} is a standard  ordinary differential equation in $V$.

Let $r>0$. Assume $u, v\in V$, $
||u||, ||v||\leq r$. Then for any $w\in V$,
\begin{equation*}
\begin{split}
\left|\int_\Omega (f_\lam(u)-f_\lam(v))wdx\right|\leq&~ \int_\Omega|f_\lam'(u+\xi v)(u-v)||w|dx\\
\leq&~(\,\mb{by $({\mb{f}_1})$}\,)\leq\int_\Omega(a_1+a_2|u+\xi v|^{p-1})|u-v||w|dx\\
\leq&~ C\int_\Omega (1+|u|^{p-1}+|v|^{p-1})|u-v||w|dx\\
\leq&~ C\(\int_\Omega (1+|u|^{p-1}+|v|^{p-1})^\frac{n}{2}dx\)^\frac{2}{n}|u-v|_{\frac{2n}{n-2}}|w|_{\frac{2n}{n-2}}\\
\leq&~ C\(\int_\Omega (1+|u|^{\frac{n}{2}(p-1)}+|v|^{\frac{n}{2}(p-1)})dx\)^\frac{2}{n}|u-v|_{\frac{2n}{n-2}}|w|_{\frac{2n}{n-2}}.
\end{split}
\end{equation*}
Observing that $p<(n+2)/(n-2)$ implies $(p-1)n/2<2n/(n-2)$,  by the Sobolev embedding $V\hookrightarrow L^{(p-1)n/2}(\Omega)$ we have
\begin{equation*}
\left|\int_\Omega (f_\lam(u)-f_\lam(v))wdx\right|\leq C(r)||u-v||||w||,
\end{equation*}
which asserts that $f:V\rightarrow V^*$ is locally Lipschitz continuous. Thus $g_\lam:V\rightarrow V$ is locally Lipschitz continuous.

 Thanks to the basic theory on abstract ODEs in Banach spaces, $(\ref{AE})$ (and hence \eqref{AE0}) has a unique  solution $u(t)$ in $V$ with initial value $u(0)=u_0$ for each $u_0\in V$. Denote by $\Phi_\lam$ the local semiflow on $V$ generated by $(\ref{AE})$, i.e., $u(t)=\Phi_\lam(t)u_0$ is the unique solution of $(\ref{AE})$ with $u(0)=u_0$.

 \br The existence results on nonclassical parabolic equations can  be also obtained by the classical  Garlerkin method; see e.g. \cite{SY}.
 \er

\subsection{Asymptotic compactness of the flow}
\Vs

In this subsection  we discuss the asymptotic compactness of $\Phi_\lam$.


        The space $V$ has an orthogonal basis $\{\varphi_j\}_{j=1}^\8$ with  $||\varphi_j||=1$ consisting of eigenvectors of $A=-\De$.
Given $m\geq 1$, denote
\begin{equation}\label{e:V1}
V_1=\mb{span}\{\varphi_1,\cdots,\varphi_m\}, \hs V_2=V_1^\bot=\text{cl}\{\mb{span}\{\varphi_j|\,\,j\geq m+1\}\},
\end{equation}
where the closure  $\text{``cl\,''}$ is taken  in $V$. Then $V=V_1\oplus V_2$.

Denote $P_m:V\rightarrow V_1$  the orthogonal projection.

\bl\label{l:4.1} Let $B$ be a bounded set in $V$. Then for any $\varepsilon>0$ there exists $m_0>0$ such that when $m>m_0$, we have
\begin{equation}\label{e£º4.4}
|(I-P_m) u|_{p+1}<\varepsilon, \Hs\forall~ u\in B.
\end{equation}
\el
\noindent{\bf Proof.} 
Let $\ve>0$ be given arbitrary. For each  $u\in V$,
since $||(I-P_m) u||\rightarrow0$ as $m\rightarrow\infty$,
 one can find a number $m_u=m_u(\ve)>0$ such that $||(I-P_m)u||<\ve$ for $m>m_u$. Hence by the Sobolev  embedding $V\hookrightarrow L^{p+1}(\Omega)$ we have
\be\label{e4.1}|(I-P_m) u|_{p+1}\leq \kappa ||(I-P_m)u||<\kappa \ve,\Hs m>m_u,\ee
where $\kappa >0$ is the embedding constant.

As the embedding $V\hookrightarrow L^{p+1}(\Omega)$ is compact,  there exist  $u_1,u_2,\cdots,u_N\in B$ such that
$B\subset\bigcup_{i=1}^N \mB_{L^{p+1}(\Omega)}\(u_i,{\varepsilon}\).$  Set $m_0=\max\{m_{u_1},\cdots,m_{u_N}\}$. Let $u\in B$. Pick a $u_j$ such that $u\in \mB_{L^{p+1}(\Omega)}\(u_j,\varepsilon\)$. Then
$$|(I-P_m)u|_{p+1}\leq |(I-P_m)(u-u_j)|_{p+1}+|(I-P_m)u_j|_{p+1}<(\kappa+1) \varepsilon,\Hs m>m_0,$$
which completes the proof of the lemma. \,$\Box$
\Vs

Let $B\subset V$. For each $u\in B$, denote
$$
T_B(u,\lam)=\sup\{t\geq 0|\,\,\,\Phi_\lam([0,t))u\subset B\}.
$$

\bl\label{l:4.2} Let $B$ be a bounded set in $V$, and $\Lam$ be a bounded interval. Then for any $\varepsilon>0$, there exist $m_0,t_0>0$  such that for all $m>m_0$ and $\lam\in\Lam$,
\begin{equation}\label{e£º4.6}
||(I-P_m)\Phi_\lam(t)u_0||<\varepsilon, \Hs\A\, u_0\in B, ~~t\in[t_0,\,T_B(u_0,\lam)).
\end{equation}

\el
\noindent{\bf Proof.} We may assume that $||u||<R$ for all $u\in B$. Let $u_0\in B$, and write  $u=u(t)=\Phi_\lam(t)u_0$. Multiplying $(\ref{AE0})$ by $u_2$ and integrating over $\Omega$, it yields
\begin{equation}\label{e:4.7}
\frac{1}{2}\frac{d}{dt}\(|u_2|^2+||u_2||^2\)+||u_2||^2=\int_\Omega f_\lam(u)u_2dx,
\end{equation}
where $u_2=(I-P_m)u.$
Using  $(\tilde{\mb{f}}_1)$ and the embedding $V\hookrightarrow L^{p+1}(\Omega)$   one can easily verify that there is a number  $C_0>0$ (independent of $\lam\in\Lam$) such that
 $$
 |f_\lam(v)|_q\leq C_0(1+||v||^{p}),\Hs v\in V,~\lam\in\Lam.
 $$

Given $\ve>0$, let  $m_0=m_0(\ve)$ be the number given by  Lemma \ref{l:4.1}. Assume $m>m_0$. Since  $u=u(t)\in B$ for $t\in[0,T_B(u_0,\lam))$,  by Lemma \ref{l:4.1} we have
\begin{equation}\label{e:4.8}\ba{ll}
\left|\int_\Omega f_\lam(u)u_2dx\right|&\leq |f_\lam(u)|_{q}|u_2|_{p+1}\leq C_0(1+||u||^p)|u_2|_{p+1}\\[1ex]
&\leq C_1|u_2|_{p+1}\leq C_1\ve,\Hs t\in[0,T_B(u_0,\lam))
\ea
\end{equation}
for all $\lam\in\Lam$.
Combining this with  $(\ref{e:4.7})$ it yields
\begin{equation*}
\frac{1}{2}\frac{d}{dt}\(|u_2|^2+||u_2||^2\)+||u_2||^2\leq C_1\ve,\Hs t\in[0,T_B(u_0,\lam)).
\end{equation*}
Observing that
$$
2||v||^2=||v||^2+||v||^2\geq {\mu_1}|v|^2+||v||^2\geq \alpha (|v|^2+||v||^2),\Hs v\in V,
$$
where $\alpha=\min\{\mu_1,1\}$, we find that
\begin{equation}\label{e:4.9}
\frac{d}{dt}\(|u_2|^2+||u_2||^2\)+\alpha(|u_2|^2+||u_2||^2)\leq 2C_1\varepsilon,\Hs ~t\in[0,T_B(u_0,\lam)).
\end{equation}
Invoking the classical Gronwall's lemma, one deduces by $(\ref{e:4.9})$ that
\begin{equation*}
\begin{split}
|u_2(t)|^2+||u_2(t)||^2\leq&~e^{-\alpha t}(|u_2(0)|^2+||u_2(0)||^2)+C_2\ve\\[1ex]
\leq&~e^{-\alpha t}\(\mu_1^{-1}+1\)||u(0)||^2+C_2\ve\\[1ex]
\leq&~ e^{-\alpha t}\(\mu_1^{-1}+1\)R^2+C_2\ve,\Hs t\in[0,T_B(u_0,\lam)),
\end{split}
\end{equation*}
where $C_2={2C_1}/{\alpha}$.
Taking $t_0={\alpha^{-1}}\log\(\(\mu_1^{-1}+1\)R^2/{\varepsilon}\)$, we immediately conclude  that
$$
|u_2(t)|^2+||u_2(t)||^2\leq \(C_2+1\)\varepsilon,\Hs t\in[t_0,T_B(u_0,\lam))
$$
which  completes the proof of the lemma. \,$\Box$

\bl\label{l:4.3} The skew-product flow $\~\Phi$ of $\{\Phi_\lam\}_{\lam\in\R}$  is asymptotically compact.
\el
\noindent{\bf Proof.} To verify the asymptotic compactness of $\~\Phi$, it suffices to check that for any $R>0$ and sequences $\lam_k\in [-R,R]$,  $u_k\in B_R:=\mB_V(0,R)$  and  $t_k\rightarrow+\8$ with $\Phi_{\lam_k}([0,t_k])u_k\subset B_R$, the sequence $v_k:=\Phi_{\lam_k}(t_k)u_k$ has a convergent subsequence.  Here (and below) $\mB_V(w,r)$ denotes the ball in $V$ centered at $w$ with radius $r$.  To this end, we only need to show that for any $\varepsilon>0$, there is a finite number of balls $\mB_V(w_i,\ve)$ $(1\leq i\leq N)$ such that
\be\label{e:6.1}\ba{ll} v_k\in \Cup_{1\leq i\leq N}\mB_V(w_i,\ve),\Hs \A\,k\geq 1.\ea\ee

We write $v_k=v_k^1+v_k^2$, where $v_k^1=P_m v_k$. As $t_k\ra+\8$, by virtue of Lemma \ref{l:4.2} there exist $m_0,k_0>0$ such that
\be\label{e:6.2}
||v_k^2||=||(I-P_m)\Phi_{\lam_k}(t_k)u_k||<{\ve}/{2},\Hs \A\,k>k_0
\ee
as long as $m>m_0$.
We fix an $m>m_0$. Then since the sequence  $v_k$ is bounded in $V$,  $v_k^1$ is a bounded sequence in $V_1=\mb{span}\{\varphi_1,\cdots,\varphi_{m}\}$.
Hence by the finite dimensionality of $V_1$ there is a finite  number of balls $\mB_V(w_i,\ve/2)$ $(1\leq i\leq N')$  such that
\be\label{e:6.3}\ba{ll} v_k^1\in \Cup_{1\leq i\leq N'}\mB_V(w_i,\ve/2),\Hs \A\,k\geq k_0.\ea\ee
Combining this with \eqref{e:6.2} one finds that
$$\ba{ll}
 v_k\in \Cup_{1\leq i\leq N'}\mB_V(w_i,\ve),\Hs \A\,k\geq k_0,\ea
 $$
 which completes the proof of \eqref{e:6.1}.  \,$\Box$

\subsection{Stability at infinity of the  flow}
\Vs

Note that the local semiflow $\Phi_\lam$ has a natural Lyapunov function
\begin{equation*}
J(u)=\frac{1}{2}||u||^2-\int_\Omega F_\lam(u)dx.
\end{equation*}
\bl\label{l:4.4} For any $c>0$, the semiflow $\Phi_\lam$  is stable in $$M_c=\{u\in V|-c\leq J(u)\leq c\}$$ at infinity in a uniform manner with respect to  $\lam$ in any bounded interval $\Lam$. Specifically, for any  $R>0$, there exists $R_1>0$ $($independent of $\lam\in \Lam$$)$ such that for any $u_0\in M_c$ with $||u_0||>R_1$, one has
\begin{equation}\label{e:4.14a}
||\Phi_\lam(t)u_0||>R,\Hs \A\,t\in\cT:=[0,T_{M_c}(u_0,\lam)),\,\,\lam\in \Lam.
\end{equation}
\el

\noindent{\bf Proof.} 
Let $u=u(t)=\Phi_\lam(t)u_0$, where $u_0\in M_c$. Multiplying $(\ref{AE0})$ by $u$ and integrating over $\Omega$, we have
\begin{equation}\label{e:4.15}
\frac{1}{2}\frac{d}{dt}\(|u|^2+||u||^2\)+||u||^2=\int_\Omega f_\lam(u)udx.
\end{equation}
By $\(\mb{f}_2\)$ there exists $\mu>2$ such that for any $\ve>0$,
\begin{equation*}
\begin{split}
\int_\Omega f_\lam(u)udx\geq&~\int_\Omega \(\mu F_\lam(u)-\varepsilon|u|^2-C_\varepsilon\)dx\\
=&~\mu\(||u||^2/2-J(u)\)-\varepsilon|u|^2-C_\varepsilon|\Omega|\\
\geq&~\(\mu/2-{\varepsilon}/{\mu_1}\)||u||^2-c\mu-C_\varepsilon|\Omega|,\Hs t\in \cT,\,\,\lam\in\Lam,
\end{split}
\end{equation*}
where $|\Omega|$ denotes the measure of $\Omega$. Therefore by $(\ref{e:4.15})$ we find that
$$
\frac{d}{dt}\(|u|^2+||u||^2\)\geq \((\mu-2)-{2\varepsilon}/{\mu_1}\)||u||^2-C'_\varepsilon,\Hs t\in\cT,\,\,\lam\in\Lam,
$$
where $C'_\varepsilon=2(c\mu+C_\varepsilon|\Omega|)$. Fix  an $\varepsilon>0$ with
$
{2\varepsilon}/{\mu_1}<(\mu-2)/2.
$
Then \begin{equation}\label{e:4.20}
\frac{d}{dt}\(|u|^2+||u||^2\)\geq \frac{1}{2}(\mu-2)||u||^2-C'_{\varepsilon},\Hs t\in \cT.
\end{equation}

Set $R_0=2\sqrt{C'_{\varepsilon}/(\mu-2)}$. By \eqref{e:4.20} one easily deduces that if $||u_0||\geq R_0$, then
$$
\frac{d}{dt}\(|u|^2+||u||^2\)\geq C'_{\varepsilon}>0,\Hs t\in \cT.
$$
It follows that
\be\label{e:6.4}|u|^2+||u||^2\geq |u_0|^2+||u_0||^2,\Hs t\in\cT.
\ee
Since $||u||^2\geq \mu_1 |u|^2$, by \eqref{e:6.4} we have
\be\label{e:4.21b}
\(\mu_1^{-1}+1\)||u||^2\geq |u_0|^2+||u_0||^2\geq ||u_0||^2,\Hs t\in\cT.
\ee
Now for any $R>0$, take $R_1=\sqrt{\mu_1+1}\,R/\sqrt{\mu_1}$. By \eqref{e:4.21b} it is easy to see  that if $||u_0||> R_1$, then  $||u||>R$ for all $t\in\cT.$
 \,$\Box$

\br The notion of stability at infinity and the techniques used here are adopted from  \cite{Li-Shi}\,$)$ and the proof of \cite[Lemma 5.5]{LLZ}, respectively.
\er
\Vs
\subsection{Proof of Theorem \ref{t:4.1}}
\Vs

Let $\Upsi$ be given as in \eqref{e:Up}. Since $\b(\lam)=f'_\lam(s)|_{s=0}$ is strictly increasing in $\lam$, by the basic knowledge on the Conley index of equilibrium solutions of evolution equations in Banach spaces (see e.g. \cite[Chap. II, Theorem 3.1]{Ryba}), we deduce  that $\Upsi$ consists of precisely all the dynamic bifurcation values of \eqref{AE} with each $\gam_k\in\Upsi$ being an essential bifurcation value; furthermore, the hypothesis (H3) in Theorem \ref{t:3.7} is fulfilled. Thus by  Theorem \ref{t:3.7}, we have

\bt\label{t:4.2} For each  $\gam_k\in \Upsilon$,  the global dynamic  bifurcation branch  $\Gamma=\Gamma(0,\gamma_k)$ of the bifurcation point $(0,\gam_k)$ $($with respect to $\{\Phi_\lam\}_{\lam\in\R}$$)$ is unbounded.
\et

 \br Nonclassical parabolic equations  arise as models to describe
physical phenomena such as non-Newtonian flow, soil mechanics and
heat conduction; see e.g. \cite{Afan}. The above theorem provides a deeper understanding to the dynamics of such equations, and is clearly of an independent interest.
  \er

We are now in a position to complete the proof of Theorem \ref{t:4.1}.\Vs

\noindent{\bf Proof of Theorem \ref{t:4.1}.} By Theorem \ref{t:4.2} the global dynamic  bifurcation branch $\Gamma=\Gamma(0,\gam_k)$  is unbounded. Set
$$
\Lambda=\{\lam|\,\,\Gamma[\lam]\ne\emp\}.
$$
Since $\Gamma$ is connected, $\Lam$ is an interval containing  $\gam_k$. By virtue of Proposition \ref{p:3.2} we know that if $\lam\in \Lam\sm\Upsi$ then $\Gamma[\lam]\ne \{0\}$. Further by Remark \ref{r:3.4} we deduce that $\Gamma[\lam]$ contains a nonempty  compact invariant set $M$ of $\Phi_\lam$ with $M\ne\{0\}$. As $\Phi_\lam$ is a gradient-like system, one easily  deduces that $M$ contains at least two distinct equilibria  of $\Phi_\lam$. Thus  $\Phi_\lam$ has at least an equilibrium $u_\lam$  in $\Gamma[\lam]$ with $u_\lam\ne 0$, which is precisely  a nontrivial  solution of \eqref{e:4.1}.

Now two cases may occur.
\vs
{\em Case} 1.) \, $\Lam$ is unbounded. In this case it is clear  that $\Lam$ contains either the interval $[\mu_k,\8)$ or the interval $(-\8,\mu_k]$, and hence
 the assertion (2) holds.
 \vs
 {\em Case }2.) \,$\Lam$ is bounded. In such a case  $\Gamma$ is  unbounded in the phase-space direction, that is, there is a bounded sequence $\lam_m\in\Lam$ such that $\sup_{u\in\Gamma[\lam_m]}||u||\ra\8$ as $m\ra\8$. For each $m$, pick a $v_m\in \Gamma[\lam_m]$ such that $||v_m||\ra\8$ as $m\ra\8$. Then we infer from Remark \ref{r:3.4} that for each $m$, there is a bounded complete trajectory $\sig_m=\sig_m(t)$ of $\Phi_{\lam_m}$ in $\Gamma[\lam_m]$ with $\sig_m(0)=v_m$.
Since $\Phi_\lam$ is a gradient system, by the basic knowledge in the theory of dynamical systems (see e.g. Hale \cite{Hale}),  both the limit sets $\a(\sig_m)$ and $\omega(\sig_m)$ of $\sig_m$ are nonempty compact invariant sets consisting of equilibrium points of $\Phi_{\lam_m}$.

 Set
 $K_m=\a(\sig_m)\cup \omega(\sig_m)$. We show that the sequence $K_m$ ($m=1,2,\cdots$) is unbounded, hence the assertion (1) holds true.

 Suppose the contrary. Then $K:=\Cup_{m\geq 1}K_m$ is a bounded set in $V$. We claim that $\sup_{u\in K}|J(u)|<\8$. Indeed, by $(\tilde{\mb{f}}_1)$ it is easy to deduce  that
\begin{equation*}
|F_\lam (s)|\leq c_1+c_2|s|^{p+1},\Hs \A\,s\in\R
\end{equation*}
for some $c_1,c_2>0$. Since $V\hookrightarrow L^{p+1}(\Omega)$, by the definition of $J$ one finds that $J$ is bounded on each bounded subset of $V$. Hence the claim holds true.

Take a positive number  $c>\sup_{u\in K}|J(u)|$. Then
\begin{equation}\label{e:4.23}
K\subset M_c=J^{-1}([-c,c]).
\end{equation}
Because  $$\max_{v\in\a(\sig_m)}J(v)\geq J(\sig_m(t))\geq \min_{v\in\omega(\sig_m)} J(v),\Hs \A\,t\in\R,$$   we see that $\sig_m$ is contained in $M_c$ for all $m$.

Take an $R>0$ so that $K\subset \mB_V(0,R)$. Then by Lemma \ref{l:4.4} there exists $R_1>2R$ such that for all $m\geq 1$ and $u_0\in M_c$ with $||u_0||>R_1$,
\begin{equation}\label{e:4.14b}
||\Phi_\lam(t)u_0||>2R,\Hs \A\,t\in [0,T_{M_c}(u_0,\lam)),\,\,\lam\in \Lam.
\end{equation}
Since $||\sig_m(0)||=||v_m||\ra\8$, \eqref{e:4.14b} implies that $$||\sig_m(t)||>2R,\Hs t\geq 0$$
for $m$ sufficiently large. Consequently $||v||\geq 2R$ for all $v\in \omega(\sig_m)$ provided that $m$ is sufficiently large, which leads to a contradiction as $\omega(\sig_m)\subset K\subset \mB_V(0,R)$ for all $m$.
The proof of Theorem $\ref{t:4.1}$ is complete.  \,$\Box$

\section{Appendix: Perturbation of A Sectorial Operator with Compact Resolvent}

Let $X$ be a Banach space, and $A$ be a  sectorial operator in $X$ with compact resolvent. For $\a\in\R$, denote $X^\a$ the fractional powers of $X$; see e.g. Henry \cite{Henry} for details.

\bp\label{p:7.1} Let $0\leq \a<1$. Assume $B:X^\a\ra X$  is a  bounded linear operator. Then  $A+B$ is a sectorial operator in $X$ with compact resolvent.
\ep
\noindent{\bf Proof.} All the argument below should be understood in the framework of complexification of spaces and operators. Since such a framework is quite standard, we omit the details.

We first observe that  $D(A)\subset X^\a=D(B)$. Let $z\in D(A)$. Then
\be\label{e:B}\ba{ll}
\|B z\|&\leq \|B\|\|z\|_\a=||B||\,||A^\a z||\\[1ex]
&\leq(\mb{by interpolation})\leq c_1||B||\,||A z||^\a||z||^{1-\a}\\[1ex]
&\leq (\mb{by Young inequality})\leq \ve ||A z||+c(\ve)||z||,
\ea
\ee
where $c(\ve)$ is a constant independent of $z$. Thus by Example (6) in Henry \cite[pp. 19]{Henry} (see also \cite[Theorem 1.3.2]{Henry}) we deduce that $A + B $ is a sectorial operator in $X$.

Note that $D(A + B )=D(A )=X^1$. To see this, by the definition of a sectorial operator it suffices to check that $A + B :D(A )\subset X\ra X$ is closed. Let $D(A )\ni z_k\ra z_0$  and $(A + B )z_k\ra w_0$ in $X$. Then
$$\ba{ll}
||A z_k-A z_m||&\leq ||(A + B )z_k-(A + B )z_m||+|| B (z_k-z_m)||\\[1ex]
&\leq (\mb{taking $\ve=\frac{1}{2}$ in \eqref{e:B}})\\[1ex]
&\leq ||(A + B )z_k-(A + B )z_m||+\frac{1}{2}||A z_k-A z_m||+c_2||z_k-z_m||.
\ea
$$
 Hence
 $$
 ||A z_k-A z_m||\leq 2||(A + B )z_k-(A + B )z_m||+2c_2||z_k-z_m||,
 $$
 which implies that $A z _k$ is a Cauchy sequence in $X$. Assume $A z_k\ra w\in X$. Then since $A $ is closed, we see that $z_0\in D(A )$ and $A z_0=w$. Using this basic fact one can easily check  that $(A + B )z_0=w_0$, which verifies  the closedness of $A + B $.

Now we show that $A + B $ has compact resolvent. For this purpose, take a number $a>0$ sufficiently large so that $\sig(A _1)\subset \mathbb{C}^+$, where $A _1=aI+A $. Then by \eqref{e:B} we find that
\be\label{e:4b}
||B A_1^{-1}||\leq \ve||AA_1^{-1}||+c(\varepsilon)||A _1^{-1}||.
\ee
It is trivial to see that $AA_1^{-1}$ is a bounded linear operator on $X$. Hence by \eqref{e:4b} one concludes that $B A_1^{-1}$ is a bounded linear operator on $X$.

Let $\lam\in\rho(A + B )$, where $\rho(A + B )$ is the resolvent set of $A + B $. We observe that
\be\label{e:4c}
\ba{ll}[\lam I-(A + B )]^{-1}=(F A _1 )^{-1},
\ea
\ee where $$F=\lam A _1^{-1}-\(AA_1^{-1}+B A_1^{-1}\)=[\lam I-(A + B )]A _1^{-1}.$$ Clearly $F:X\ra X$ is bounded.
Since both operators $A _1^{-1}:X\ra D(A )$ and $\lam I-(A + B ): D(A + B )=D(A )\ra X$ are one-one mappings, we deduce that $F$ is a one-one mapping on $X$. The classical Banach inverse theorem then asserts that $F^{-1}$ is a bounded linear operator on $X$. Thus by \eqref{e:4c} we deduce  that
\begin{equation}\label{e:4d}
[\lam I-(A + B )]^{-1}=A _1^{-1}F^{-1}.
\end{equation}
Since  $A _1^{-1}:X\ra X^1$ is compact, by \eqref{e:4d} we immediately conclude that $[\lam I-(A + B )]^{-1}$ is compact. $\Box$




\Vs

\end{document}